\documentclass[a4paper,reqno,makeidx]{amsart}
\usepackage{a4wide}
\RequirePackage{ifthen}
\newboolean{omarfont}
{%
  \setboolean{omarfont}{true}%
}{%
  \setboolean{omarfont}{false}%
}
  \RequirePackage{ifthen}
  \provideboolean{useutopia}
  \setboolean{useutopia}{false}%
  \ifthenelse{\boolean{useutopia}}{%
    \RequirePackage{ifthen}
    \RequirePackage{iftex}
    \provideboolean{usemathrsfs}%
    \provideboolean{useutopia}
    \setboolean{useutopia}{false}%
    \ifXeTeX
      \RequirePackage[british]{babel}
      \RequirePackage{xltxtra}%
      \ifthenelse{\boolean{useutopia}}{
        \RequirePackage[utopia,euro]{mathdesign}
        \RequirePackage[OMLmathbf,OMLmathsfit]{isomath}
      }{
        \RequirePackage{amssymb}    %
        \RequirePackage{bbold}    %
      }
        \RequirePackage{mathrsfs} %
        \setboolean{usemathrsfs}{true}%
      \else
      \RequirePackage[utf8]{inputenc}
      \ifthenelse{\boolean{useutopia}}{%
        \RequirePackage[utopia,euro]{mathdesign}%
        \RequirePackage[OMLmathbf,OMLmathsfit]{isomath}
        \DeclareSymbolFont{usualmathcal}{OMS}{cmsy}{m}{n}
        \DeclareSymbolFontAlphabet{\mathcalbf}{usualmathcal}
        \RequirePackage{stmaryrd} %
        \providecommand{\diracdelta}[1][]{\ensuremath{\deltaup_{#1}}}

        \providecommand{\measure}[1]{\ensuremath{\mathcalbf{\uppercase{#1}}}}
        \providecommand{\olnum}[1]{\ensuremath{\mathsfit{#1}}}
        \providecommand{\numsca}[1]{\olnum{#1}}
        \providecommand{\numvec}[1]{\ensuremath{\mathsfbfit{#1}}}%
      }{%
        \RequirePackage{amssymb}    %
        \RequirePackage{mathrsfs} %
        \RequirePackage{bbold}    %
        \RequirePackage{stmaryrd} %
        \setboolean{usemathrsfs}{true}%
        \providecommand{\mathcalbf}{\mathcal}
      }%
      \RequirePackage[british]{babel}
    \fi
    \RequirePackage{mathtools}%
    \RequirePackage{nicefrac}
    \RequirePackage{amsthm}
    \RequirePackage{enumerate}%
    \RequirePackage{xspace}%
    \RequirePackage{verbatim}%
    \RequirePackage{fancyvrb}
    \RequirePackage{hyperref}%
    \RequirePackage{listings}%
    \RequirePackage{xcolor}
    \RequirePackage{booktabs}%
    \RequirePackage{tikz}%
    \usetikzlibrary{calc}%
    \provideboolean{usebeamer}%
    \provideboolean{isthesis}%
    \provideboolean{isamsltex}%
    \provideboolean{issiamltex}%
  }{%
    \RequirePackage{ifthen}
    \RequirePackage{iftex}
    \provideboolean{usemathrsfs}%
    \provideboolean{useutopia}
    \setboolean{useutopia}{false}%
    \ifXeTeX
      \RequirePackage[british]{babel}
      \RequirePackage{xltxtra}%
      \ifthenelse{\boolean{useutopia}}{
        \RequirePackage[utopia,euro]{mathdesign}
        \RequirePackage[OMLmathbf,OMLmathsfit]{isomath}
      }{
        \RequirePackage{amssymb}    %
        \RequirePackage{bbold}    %
      }
        \RequirePackage{mathrsfs} %
        \setboolean{usemathrsfs}{true}%
      \else
      \RequirePackage[utf8]{inputenc}
      \ifthenelse{\boolean{useutopia}}{%
        \RequirePackage[utopia,euro]{mathdesign}%
        \RequirePackage[OMLmathbf,OMLmathsfit]{isomath}
        \DeclareSymbolFont{usualmathcal}{OMS}{cmsy}{m}{n}
        \DeclareSymbolFontAlphabet{\mathcalbf}{usualmathcal}
        \RequirePackage{stmaryrd} %
        \providecommand{\diracdelta}[1][]{\ensuremath{\deltaup_{#1}}}

        \providecommand{\measure}[1]{\ensuremath{\mathcalbf{\uppercase{#1}}}}
        \providecommand{\olnum}[1]{\ensuremath{\mathsfit{#1}}}
        \providecommand{\numsca}[1]{\olnum{#1}}
        \providecommand{\numvec}[1]{\ensuremath{\mathsfbfit{#1}}}%
      }{%
        \RequirePackage{amssymb}    %
        \RequirePackage{mathrsfs} %
        \RequirePackage{bbold}    %
        \RequirePackage{stmaryrd} %
        \setboolean{usemathrsfs}{true}%
        \providecommand{\mathcalbf}{\mathcal}
      }%
      \RequirePackage[british]{babel}
    \fi
    \RequirePackage{mathtools}%
    \RequirePackage{nicefrac}
    \RequirePackage{amsthm}
    \RequirePackage{enumerate}%
    \RequirePackage{xspace}%
    \RequirePackage{verbatim}%
    \RequirePackage{fancyvrb}
    \RequirePackage{hyperref}%
    \RequirePackage{listings}%
    \RequirePackage{xcolor}
    \RequirePackage{booktabs}%
    \RequirePackage{tikz}%
    \usetikzlibrary{calc}%
    \provideboolean{usebeamer}%
    \provideboolean{isthesis}%
    \provideboolean{isamsltex}%
    \provideboolean{issiamltex}%
  }
  \colorlet{a}{magenta}
  \colorlet{b}{green!70!blue}
  \colorlet{c}{yellow!50!black}
  \colorlet{d}{cyan}
  \colorlet{e}{red}
  \colorlet{f}{blue}
  \colorlet{g}{white}
  \colorlet{h}{black!50}
  \colorlet{i}{black}
  \colorlet{j}{black!75}

  \providecommand{\linkedurl}[1]{\url{#1}}

  \providecommand{\Ignore}[1]{}
  \providecommand{\ignore}[1]{}
  \providecommand{\freeze}[1]{}%
  \providecommand{\crossout}[1]{{\textcolor{i!20}{#1}}}
  \providecommand{\highlight}[1]{{\color{a}#1}}

  \providecommand{\memo}[1]{
    \ensuremath{
      \framebox{\tiny\textbf{\kern-2pt\textsf{#1}}\kern-2pt}
    }
    \xspace}

  \RequirePackage{alphalph}
  \provideboolean{shownotes}
  \setboolean{shownotes}{true}%
  \newcounter{margnote}[page]
  \providecommand{\mgcolor}{a}
  \providecommand{\mgcolorset}[1]{\renewcommand{\mgcolor}{\alphalph{#1}}}
  \providecommand{\mgcolorsetbycounter}[1]{\renewcommand{\mgcolor}{\alph{#1}}}
  \providecommand{\mgcolormake}{\mgcolorsetbycounter{margnote}}
  \providecommand{\mgcolorstepby}[1]{
    \setcounter{tmpcounter}{\value{margnote}}%
    \addtocounter{tmpcounter}{#1}%
    \mgcolorsetbycounter{tmpcounter}%
  }%
  \providecommand{\margnotecolor}{
    \ifthenelse{\value{margnote}<7}{%
      \mgcolormake%
    }{%
      \ifthenelse{\value{margnote}=7}{\mgcolorset{10}}{%
        \ifthenelse{\value{margnote}<11}{\mgcolormake}{%
          \ifthenelse{\value{margnote}<17}{\mgcolorstepby{-10}}{%
              \mgcolorset{10}%
          }%
        }%
      }%
    }%
  }%
  \providecommand{\margnotemark}{{\colorbox{\mgcolor}{\tiny\color{g}\upshape\texttt{\arabic{margnote}}}}}
  \providecommand{\margnote}[2][]{%
    \ifthenelse{%
      \boolean{shownotes}%
    }{%
      \stepcounter{margnote}%
      \margnotecolor%
      \margnotemark%
      \marginpar{%
        \color{\mgcolor}%
        \texttt{%
          \begin{minipage}{2cm}%
            \raggedright\tiny%
            \margnotemark%
                {\ifx|#1|{}\else{#1:\ }\fi}%
                #2%
          \end{minipage}%
        }%
      }%
    }{%
    }%
  }%
  \providecommand{\mathnote}[2][]{%
    \ifthenelse{%
      \boolean{shownotes}%
    }{%
      \stepcounter{margnote}%
      \margnotecolor%
      \text{%
        \colorbox{\mgcolor}{%
          \color{g}%
          \texttt{%
            \tiny%
                [\arabic{margnote}]%
                \ifx|#1|{}\else{#1:}\fi%
                #2%
          }%
        }%
      }%
    }{%
    }%
  }%
  
  \providecommand{\textnote}[2][]{%
    \ifthenelse{%
      \boolean{shownotes}%
    }{%
      \stepcounter{margnote}%
      \margnotecolor%
      \ \\
      \text{%
        \colorbox{\mgcolor}{%
          \begin{minipage}{.9\textwidth}
          \color{g}%
          \texttt{%
            \tiny%
            [\arabic{margnote}]%
            \ifx|#1|{}\else{#1: }\fi%
            #2%
          }%
          \end{minipage}
        }%
      }%
    }{%
    }%
  }%

  \providecommand{\todo}[1]{\ifthenelse{\boolean{shownotes}}{\margnote[To do.]{#1}}{}}
  \providecommand{\Todo}[1]{
    \ifthenelse{\boolean{shownotes}}{
      \begin{center}
      \begin{tikzpicture}
       \node[fill=a!17]{
         \begin{minipage}{\textwidth}
           \texttt{To do:}
           \\
           \texttt{\bfseries{\small #1}}
         \end{minipage}
       };
      \end{tikzpicture}
      \end{center}
    }{}}
  \provideboolean{showrevisions}
  \setboolean{showrevisions}{true}
  \newcommand{\revisionsheader}{***\newline\Warning{the following part is under revision}}
  \newcommand{\revisionsfooter}{***\newline\Warning{end of part under revision}}

  \providecommand{\Warning}[1]{    
    \begin{tikzpicture}
      \node[fill=a!27]{
        \begin{minipage}{\textwidth}
          \texttt{\bfseries{\small Warning: #1}}
        \end{minipage}
      };
    \end{tikzpicture}
  }

  \provideboolean{showcomments}
  \providecommand{\margincomment}[1]{
  \ifthenelse{\boolean{showcomments}}{\marginpar{\tiny #1}}{}
  }
  \provideboolean{showchanges}
  \setboolean{showchanges}{false}
  \providecommand{\changes}[1]{
    \ifthenelse{\boolean{showchanges}}{{\highlight{#1}}}{#1}
  }
  \providecommand{\changefromto}[3][replace with]{%
    \ifthenelse{\boolean{showchanges}}{%
      {\crossout{#2}\margnote{#1}}{\highlight{#3}}}{%
      #3\xspace}%
  }
  \providecommand{\ChangePar}[2]{
    \ifthenelse{\boolean{showchanges}}
    {{\par\textcolor{i!20}{#1}}{\par\textcolor{a}{#2}}}
    {\par #2}
  }
  \providecommand{\InsertPar}[1]{
    \ifthenelse{\boolean{showchanges}}
    {{\par$\mapsto$ \textcolor{blue}{#1}}}
    {\par #1}
  }

  \providecommand{\mathscript}
  	   {\mathscr}

   \providecommand{\cB}{\ensuremath{\mathscript B}\xspace}
   \providecommand{\cC}{\ensuremath{\mathscript C}\xspace}

   \providecommand{\cF}{\ensuremath{\mathscript F}\xspace}
   
   \providecommand{\cH}{\ensuremath{\mathscript H}\xspace}

   \providecommand{\cV}{\ensuremath{\mathscript V}\xspace}
   
   \providecommand{\cX}{\ensuremath{\mathscript X}\xspace}

   \providecommand{\bbbold}{\mathbb}

   \providecommand{\rN}{\ensuremath{\bbbold N}\xspace}
   
   \providecommand{\rP}{\ensuremath{\bbbold P}\xspace}
   
   \providecommand{\rR}{\ensuremath{\bbbold R}\xspace}

  \providecommand{\Ae}[1][]{\ensuremath{\ifx|#1|{\ }\else{\:#1\text{-}}\fi\text{almost everywhere }}\xspace}
  \providecommand{\Aa}[1][]{\ensuremath{\text{ for }\ifx|#1|{}\else{\:#1\text{-}}\fi\text{almost all }}}
  \providecommand{\as}[1][]{\ensuremath{\ifx|#1|{\ }\else{#1\text{-}}\fi\text{almost surely}}\xspace}
   \providecommand{\naturals}{\rN\xspace}
   
   \providecommand{\NO}{\ensuremath{\naturals_0}}

   \providecommand{\reals}{\rR}
   \providecommand{\setvecs}[2]{{#1}^{#2}}
   \providecommand{\setmats}[3]{\setvecs{#1}{#2\times{#3}}}
   \providecommand{\psetvecs}[2]{\powqp{#2}{#1}}
   \providecommand{\psetmats}[3]{\psetvecs{#1}{#2\times{#3}}}
   \providecommand{\R}[1]{\reals^{#1}}
   
   \providecommand{\fieldmats}[3][F]{\csname#1\endcsname{#2\times#3}}
   \providecommand{\realmats}[2]{\fieldmats[R]{#1}{#2}}

   \providecommand{\one}{\ensuremath{\bbbold 1}\xspace}
   
   \providecommand{\iverson}[1]{\one_{\qb{#1}}}

   \providecommand{\diracdelta}[1][]{\ensuremath{{\mathrm{\delta}}\ifx|#1|{}\else_{#1}\fi}}

   \providecommand{\pic}{\ensuremath{\mathrm\pi}}

   \providecommand{\Tolto}{\smallsetminus}
   \providecommand{\take}{\Tolto}
   
   \providecommand{\closure}[1]{\overline{#1}}
   \providecommand{\inner}{\cdot}

   \providecommand{\frobinner}{\!:\!}

   \providecommand{\Y}{\ensuremath{\varUpsilon}\xspace}

   \providecommand{\W}{\ensuremath{\varOmega}\xspace}
   \providecommand{\qgroup}[1]{{#1}}%
   
   \providecommand{\qp}[1]{\ensuremath{\left({#1}\right)}}
   \providecommand{\qpreg}[1]{\ensuremath{(#1)}}
   \providecommand{\qpbig}[1]{\ensuremath{\big(#1\big)}}
   \providecommand{\qpBig}[1]{\ensuremath{\Big(#1\Big)}}
   \providecommand{\qpbigg}[1]{\ensuremath{\bigg(\!#1\!\bigg)}}
   \providecommand{\qpBigg}[1]{\ensuremath{\Bigg(\!#1\!\Bigg)}}
   \providecommand{\qb}[1]{\ensuremath{\left[{#1}\right]}}

   \providecommand{\qc}[1]{\ensuremath{\left\{{#1}\right\}}}

   \providecommand{\qa}[1]{\ensuremath{\left\langle{#1}\right\rangle}}

   \providecommand{\opinter}[2]{\ensuremath{\left(#1,#2\right)}\xspace}
   \providecommand{\clinter}[2]{\ensuremath{\left[#1,#2\right]}\xspace}

   \providecommand{\compowqp}[2]{\ensuremath{\qp{\!#2\!\!}^{\kern -.4em #1}\!}}
   
   \providecommand{\powqpreg}[2]{\ensuremath{%
       \qpreg{#2}^{\kern 0em\lower .1ex\hbox{\scriptsize $#1$}}\kern-.3em}}
   \providecommand{\powqpbig}[2]{\ensuremath{%
       \qpbig{#2}^{\kern -.2em\lower .3ex\hbox{\scriptsize $#1$}}\kern-.3em}}
   \providecommand{\powqpBig}[2]{\ensuremath{%
       \qpBig{#2}^{\kern -.2em\lower .3ex\hbox{\scriptsize $#1$}}\kern-.3em}}
   \providecommand{\powqpbigg}[2]{\ensuremath{%
       \qpbigg{#2}^{\kern -.2em\lower .3ex\hbox{\scriptsize $#1$}}\kern-.3em}}
   \providecommand{\powqpBigg}[2]{\ensuremath{%
       \qpBigg{#2}^{\kern -.2em\lower .3ex\hbox{\scriptsize $#1$}}}}
   \providecommand{\powp}[3][]{#3\ifx|#1|^{#2}\else{#1}^{#2}\fi}%
   \providecommand{\pow}[2][]{\ifx|#1|\operatorname{pow}^{#2}\else\powp{#2}{#1}\fi}%

   \providecommand{\powqp}[3][]{\powp[#1]{#2}{\qp{#3}}}%

   \providecommand{\powabs}[2]{\powp{#1}{\abs{#2}}}
   
   \providecommand{\pownorm}[2]{\powp{#1}{\norm{#2}}}

   \providecommand{\norm}[2][]{\ifx|#1|\left|\else\csname#1\endcsname|\fi#2\ifx|#1|\right|\else\csname#1\endcsname|\fi}

   \providecommand{\abs}[2][]{\ensuremath{\ifx|#1|{\left|#2\right|}\else{\csname#1\endcsname|{#2}\csname#1\endcsname|}\fi}}

   \providecommand{\Norm}[1]{\ensuremath{\left\|#1\right\|}}
   \providecommand{\ltwop}[2]{\ensuremath{\qa{#1,#2}}}

   \providecommand{\duality}[3][]{\ensuremath{#1\langle #2\,#1\vert\,#3#1\rangle}}
   \providecommand{\average}[2][]{{\qa{#2}\ifx|#1|\else_{#1}\fi}}

   \providecommand{\ensemble}[2]{\ensuremath{\left\{ #1:\;#2 \right\}}}
   \providecommand{\setofsuch}{\ensemble}%

   \providecommand{\setof}[1]{{\qc{#1}}}

   \providecommand{\conditionalto}[1]{{\left|{#1}\right.}}

  \providecommand{\measure}[1]{\ensuremath{\mathcal{\uppercase{#1}}}}
  
  \providecommand{\probmeasure}[2][]{{\measure{#2}}_{#1}}
  \providecommand{\Prob}{}
  \renewcommand{\Prob}[1][]{\probmeasure[{#1}]{p}}

  \providecommand{\randvars}[1][\Prob]{\operatorname{RV}\ifx|#1|{}\else{(#1)}\fi}
  \providecommand{\discrandvars}[1][\Prob]{\operatorname{DRV}\ifx|#1|{}\else{({#1)}\fi}} 
  \providecommand{\contrandvars}[1][\Prob]{\ensuremath{\operatorname{CDRV}\ifx|#1|{}\else(#1)\fi}} 
   \def\env@matrix{\hskip -\arraycolsep
    \let\@ifnextchar\new@ifnextchar
    \array{*\c@MaxMatrixCols c}}
   \makeatletter
   \renewcommand*\env@matrix[1][c]{\hskip -\arraycolsep
     \let\@ifnextchar\new@ifnextchar
     \array{*\c@MaxMatrixCols #1}}
   \makeatother
   \providecommand{\irow}[2]{#1_{#2}}%
   \providecommand{\icol}[2]{#1^{#2}}%

   \providecommand{\ijrowcol}[3]{\icol{\irow{#1}{#2}}{#3}}
   \providecommand{\entry}[1]{\qb{#1}}
   \providecommand{\vecentry}[2]{\irow{#1}{#2}}

   \providecommand{\vecof}[1]{\qp{#1}}
   
   \providecommand{\rowof}[1]{\qb{#1}}

   \providecommand{\getentryi}[2]{\irow{\entry{#1}}{#2}}
   
   \providecommand{\getvecentry}[2]{\getentryi{\vec #1}{#2}}

   \providecommand{\vecidotsfromto}[3]{\vecof{\listidotsfromto{#1}{#2}{#3}}}

   \providecommand{\getrow}[2]{\irow{\entry{#1}}{#2}}

   \providecommand{\dismatof}[2][r]{\begin{bmatrix}[#1]#2\end{bmatrix}}

   \providecommand{\getentryij}[3]{\ijrowcol{\entry{#1}}{#2}{#3}}
   \providecommand{\matentry}[3]{\ijrowcol{#1}{#2}{#3}}

   \providecommand{\block}[5]{\ijrowcol{#1}{\ifx#2#3{\rowof{#2}}\else\rowof{{#2}\dotsc{#3}}\fi}{\ifx#4#5{\rowof{#4}}\else\rowof{{#4}\dotsc{#5}}\fi}}
   \providecommand{\colblock}[3]{\getvecentry{#1}{\ifx#2#3{#2}\else\fromto{#2}{#3}\fi}}

   \providecommand{\matijnm}[6][1]{\ensuremath{\ijrowcol{\smash{\qb{\ijrowcol {#2}{#3}{#4}}}}{\rangefromto {#3}{#1}{#5}}{\rangefromto {#4}{#1}{#6}}}}

   \providecommand{\dismatskeldots}[4]{
     \dismatof[c]{
       #1&\dotsc&#3
       \\
       \vdots & \ddots &\vdots
       \\
       #2&\dotsc&#4
     }
   }
   \providecommand{\dismatcommfromtofromto}[5]{
     \dismatskeldots{#1#2#4}{#1#3#4}{#1#2#5}{#1#3#5}
   }
   \providecommand{\dismatcustfromtofromto}[6][matentry]{
     \dismatcommfromtofromto{\csname#1\endcsname{#2}}#3#4#5#6
   }
   \providecommand{\dismatcustfromtofromto}[6][matentry]{
     \dismatskeldots{%
       \csname#1\endcsname{#2}{#3}{#4}%
     }{%
       \csname#1\endcsname{#2}{#3}{#6}%
     }{%
       \csname#1\endcsname{#2}{#5}{#4}%
     }{%
       \csname#1\endcsname{#2}{#5}{#6}%
     }%
   }%
   \providecommand{\dismatcustfromtofromto}[6][matentry]{
     \dismatof{
       \csname#1\endcsname{#2}{#3}{#4}&\dotsc&\csname#1\endcsname{#2}{#3}{#6}
       \\
       \vdots & \ddots &\vdots
       \\
       \csname#1\endcsname{#2}{#5}{#4}&\dotsc&\csname#1\endcsname{#2}{#5}{#6}
     }
   }
   \providecommand{\dismatcustnm}[4][matentry]{\dismatcustfromtofromto[#1]{#2}1{#3}1{#4}}
   
   \providecommand{\squmatcustmd}[3][matentry]{\dismatcustnm[#1]{#2}{#3}{#3}}
   \providecommand{\squmatmd}[2]{\squmatcustmd[matentry]{#1}{#2}}
   
   \providecommand{\dissysaxbdotsnm}[5]{\begin{matrix}[r]%
       \matentry{#1}11\vecentry{#2}1&+\dotsb&+\matentry{#1}1{#5}\vecentry{#2}{#5}
       &
       =
       \ifx|#3|0\else{\vecentry {#3}1}\fi
       \\
       \dotsb
       \\
       \matentry{#1}{#4}1\vecentry{#2}1&+\dotsb&+\matentry{#1}{#4}{#5}\vecentry{#2}{#5}
       &
       =
       \ifx|#3|0\else{\vecentry {#3}{#4}}\fi
   \end{matrix}}
   \providecommand{\seqof}[1]{\qp{#1}}%
   \providecommand{\seq}[1]{\seqof{#1}}%
   \providecommand{\seqs}[2]{\seqof{#1}_{#2}}
   \providecommand{\sets}[2]{\setof{#1}_{#2}}%
   \providecommand{\seqi}[3][]{\seqs{#2_{#3}}{\ifx|#1|{#3}\else{{#3}\in{#1}}\fi}}%

   \providecommand{\seti}[3][]{\sets{#2_{#3}}{\ifx|#1|_{#3}\else_{{#3}\in{#1}}\fi}}%
   \providecommand{\sequ}[3][]{\seqs{#2^{#3}}{\ifx|#1|{#3}\else{{#3}\in{#1}}\fi}}%
   \providecommand{\setu}[3][]{\sets{#2^{#3}}{\ifx|#1|{#3}\else{{#3}\in{#1}}\fi}}%
   \providecommand{\limofat}[3][]{\ensuremath{\lim_{\ifx|#1|{}\else{#1\ni}\fi#3}{#2}}}

   \providecommand{\listitwo}[1]{\ensuremath{#1_1,#1_2}}
   \providecommand{\listidotsfromto}[3]{\ensuremath{#1_{#2},\dotsc,#1_{#3}}}

   \providecommand{\seqidotsfromto}[3]{\seq{\listidotsfromto{#1}{#2}{#3}}}

   \providecommand{\sumifromto}[3]{\ensuremath{\sum_{#1=#2}^{#3}}}

   \providecommand{\jump}[2][]{\ensuremath{\left\llbracket #2\right\rrbracket\ifx|#1|{}\else_{#1}\fi}}

   \providecommand{\fromto}[2]{\ensuremath{\setof{#1\dotsc#2}}}%

   \providecommand{\integerbetween}[2]{\ensuremath{={#1},\dotsc,{#2}}}
   \providecommand{\rangefromto}[3]{\ensuremath{#1\integerbetween{#2}{#3}}}

   \providecommand{\d}{}
   \renewcommand{\d}[1][]{\ensuremath{\operatorname{d}\!\ifx|#1|\else{_{#1}}\fi}}
   
   \providecommand{\ds}[1][]{\d{\measure S}}
   \providecommand{\D}[1][]{\ensuremath{\operatorname{D}\!\ifx|#1|\else{_{#1}}\fi}}

  \providecommand{\registered}%
  {\ensuremath{^\text{\textregistered}}}

  \providecommand{\tand}{\ensuremath{\text{ and }}}

  \providecommand{\constant}[1]{\ensuremath{C_{#1}}}
  \providecommand{\constext}[2][]{\constant{\text{#2}{\ifx|#1|{}\else{,\ensuremath{#1}}\fi}}}            %
  \providecommand{\constref}[2][]{\ensuremath{\constant{\textup{\ref{#2}{\ifx|#1|{}\else{,\ensuremath{#1}}\fi}}}}}
  \providecommand{\constdef}[2][]{\label{#2}\ensuremath{\constant{\textup{\ref{#2}{\ifx|#1|{}\else{,\ensuremath{#1}}\fi}}}}}

  \providecommand{\funkref}[3][]{\ensuremath{{#3}_{\textup{\ref{#2}{\ifx|#1|{}\else{,\ensuremath{#1}}\fi}}}}}

  \providecommand{\Cof}{\operatorname{Cof}}

  \providecommand{\Diag}{\operatorname{Diag}}
  
  \providecommand{\diam}{\operatorname{diam}}
  \providecommand{\dist}{\operatorname{dist}}
  
  \renewcommand{\div}[1][]{\nabla\ifx|#1|{}\else\kern-2pt_{#1}\fi\kern-2pt\cdot}
  \providecommand{\divof}[2][]{\div[#1]\ifx|#2|{}\else\qb{#2}\fi}

  \providecommand{\inverse}[2][]{\powp[#1]{-1}{#2}}
  \providecommand{\inverseqp}[1]{\inverse{\qp{#1}}}
  \providecommand{\inverseof}[1]{\inverseqp{#1}}
  
  \providecommand{\inversemat}[1]{\inverse{\mat{#1}}}

  \providecommand{\fracl}[3][]{\ifx|#1|\nicefrac{#2}{#3}\else{#2}#1/{#3}\fi}
  \providecommand{\qpfracl}[3][]{\qp{\ifx|#1|\fracl{#2}{#3}\else{#2}#1/{#3}\fi}}
  \providecommand{\qpfrac}[3][]{\qp{\ifx|#1|\frac{#2}{#3}\else{#2}#1/{#3}\fi}}
  \providecommand{\absfracl}[3][]{\abs{\ifx|#1|\fracl{#2}{#3}\else{#2}#1/{#3}\fi}}
  \providecommand{\absfrac}[3][]{\abs{\ifx|#1|\frac{#2}{#3}\else{#2}#1/{#3}\fi}}
  \providecommand{\fraclff}[3][]{\ifx|#1|{#2}/{#3}\else{#2}#1/{#3}\fi}

  \providecommand{\eye}[1][]{\boldsymbol{\mathrm I}\ifx|#1|{}\else_{#1}\fi}%
  \providecommand{\numeye}[1][]{\boldsymbol{\mathsf{I}}\ifx|#1|{}\else_{#1}\fi}%
  \providecommand{\Eye}[1]{
    \begin{bmatrix}
    \ifthenelse{#1>1}{
      \ifthenelse{#1>2}{
        \ifthenelse{#1>3}{
          \ifthenelse{#1>4}{
            1&\zeroentry&\dotso&\zeroentry
            \\
            \zeroentry&1&\dotso&\zeroentry
            \\
            \vdots&\vdots&\ddots&\vdots
            \\
            \zeroentry&\zeroentry&\dotso&1
          }{        
            1&\zeroentry&\zeroentry&\zeroentry
            \\
            \zeroentry&1&\zeroentry&\zeroentry
            \\
            \zeroentry&\zeroentry&1&\zeroentry
            \\
            \zeroentry&\zeroentry&\zeroentry&1
          }
        }{
          1&\zeroentry&\zeroentry
          \\
          \zeroentry&1&\zeroentry
          \\
          \zeroentry&\zeroentry&1
        }
      }{
        1&\zeroentry
        \\
        \zeroentry&1
      }
    }{
      1
    }
    \end{bmatrix}
  }

  \providecommand{\lebmeas}[1][]{\operatorname{l}^{#1}}     %
  \providecommand{\lebmeasof}[2][]{\ifx|#1|\left|#1\right|\else\lebmeas[#1]\qp{#2}\fi}         %
  \providecommand{\dash}[1][']{\ifthenelse{\equal{#1}{'}\OR\equal{#1}{''}}{#1}{^{(#1)}}}
  
  \providecommand{\pdfrac}[2][]{\ensuremath{\frac{\partial\ifx|#1|\phantom{#2}\else{#1}\fi}{\partial{#2}}}} %
  \providecommand{\pd}[2][]{\ensuremath{\partial_{#2}}{\ifx|#1|{}\else{\qb{#1}}\fi}} %

  \renewcommand{\Im}{\operatorname{im}}                 %
  \renewcommand{\Re}{\operatorname{re}}                 %
  \providecommand{\imaginpart}[1][]{\Im{\ifx|#1|{}\else\qp{#1}\fi}} %
  \providecommand{\realpart}[1][]{\Re{\ifx|#1|{}\else\qp{#1}\fi}} %
  \providecommand{\trace}{\operatorname{tra}}             %

  \providecommand{\transpose}{\intercal}%

  \providecommand{\transposeof}[1]{\ensuremath{\qp{#1}^\transpose}}
  \providecommand{\transposed}{{}^\transpose}
  
  \providecommand{\Transpose}[1]{\ensuremath{{#1}^{\transpose}}}

  \providecommand{\transnumvec}[1]{\Transpose{\numvec{#1}}}
  \providecommand{\Transposemat}[1]{\Transpose{\mat{#1}}}

  \providecommand{\transposemat}[1]{\Transposemat{#1}}
  \providecommand{\Transinverse}[1]{\Transpose{{{#1}^{-}}}} %
  \providecommand{\transinversemat}[1]{\Transinverse{\mat{#1}}}

  \providecommand{\orthogonalto}[1][]{\ensuremath{\perp\ifx|#1|{}\else{_{#1}}\fi}}
  
  \providecommand{\rowof}[1]{\ensuremath{\vecof{#1}}}

  \providecommand{\colvec}[1]{\ensuremath{\vecof{#1}}}

  \providecommand{\colvecthree}[3]{\ensuremath{\colvec{#1,#2,#3}}}
  \providecommand{\colvecfour}[4]{\ensuremath{\colvec{#1,#2,#3,#4}}}

  \providecommand{\vecthree}{\colvecthree}
  \providecommand{\vecfour}{\colvecfour}

  \providecommand{\disrowof}[1]{\rowof{\begin{matrix}[r]#1\end{matrix}}}
  \providecommand{\disrowvec}[1]{\disrowof{#1}}
  \providecommand{\zeroentry}{\phantom0}%

  \providecommand{\disrowvecfour}[4]{\ensuremath{\disrowvec{#1 &#2 &#3 &#4}}}

  \providecommand{\rowvecdotsfromto}[3]{\rowof{\icol{#1}{#2},\dotsc,\icol{#1}{#3}}}

  \providecommand{\inta}[1]{\qa{#1}}
  \providecommand{\inton}[2]{\inta{#1}_{#2}}
  \providecommand{\normalsymbol}{\operatorname{n}}
  \providecommand{\normal}[1][]{\normalsymbol\ifx|#1|\else_{#1}\fi}%
  \providecommand{\normalto}[1]{\ensuremath{\normal[#1]}}
  \providecommand{\normalder}[1][]{\ensuremath{\normal\ifx|#1|\else\qp{#1}\fi{\inner\grad}}}
  \providecommand{\normalderto}[2][]{\ensuremath{\normalto{#2}\ifx|#1|\else\qp{#1}\fi{\inner\grad}}}

  \providecommand{\tangentialsymbol}{\operatorname{t}}
  \providecommand{\tangentialto}[2][]{\tangentialsymbol\ifx|#1|\else^{#1}\fi\ifx|#2|\else_{#2}\fi}

  \providecommand{\intersected}{\ensuremath{\cap}}
  
  \providecommand{\meet}{\intersected}

  \renewcommand{\vec}[1]{\ensuremath{\boldsymbol{#1}}}
  \providecommand{\hatvec}[1]{\hat{\vec{#1}}}
  \providecommand{\hatmat}[1]{\hat{\mat{#1}}}

  \providecommand{\geomat}[1]{\vec{#1}}

  \providecommand{\mat}[1]{\geomat{#1}} %
  \providecommand{\zeromat}{\ensuremath{\boldsymbol{\mathrm O}}}           %
  \providecommand{\olnum}[1]{\ensuremath{\mathsf{#1}}} %
  \providecommand{\numsca}[1]{\olnum{#1}}
  \providecommand{\numvec}[1]{{\vec{\numsca{#1}}}}
  \providecommand{\numvecentry}[2]{\irow{\numsca{#1}}{#2}}
  
  \providecommand{\nummat}[1]{\numvec{\MakeUppercase{#1}}}

  \providecommand{\transnumvec}[1]{\Transpose{\numvec{#1}}}
  \providecommand{\Prob}[1][]{\ensuremath{\operatorname{Prob}\ifx|#1|{}\else_{#1}\fi}}

  \providecommand{\pdf}[2][]{\ensuremath{\operatorname{pdf}_{#2\ifx|#1|{}\else{\conditionalto{#1}}\fi}}\xspace}

  \providecommand{\expectation}{\ensuremath{\operatorname{E}}}
  \providecommand{\EX}[1][]{\ensuremath{\expectation\ifx|#1|{}\else_{#1}\fi}}

  \providecommand{\gausskernel}[3][x]{%
    \ensuremath{
      \exp\frac{-\if#20{#1}\else(#1-\mu)\fi^2}{%
        2\if#31{}\else\powp2{#3}\fi}%
    }%
  }
  \providecommand{\gaussdistribution}[3][x]{%
    \ensuremath{\frac1{\sqrt{2\pic}\if#31{}\else#3\fi}%
      \gausskernel[#1]{#2}{#3}
    }%
  }%

  \providecommand{\boundary}{\partial}

  \providecommand{\SPD}{\operatorname{SPD}}
  \providecommand{\spdmats}[2][F]{\SPD(\csname#1\endcsname{#2})}
   \providecommand{\Continuous}{\ensuremath{\operatorname C}\xspace}%
   \providecommand{\Hspace}{\ensuremath{\operatorname H}\xspace}
   \providecommand{\Lebesgue}{\ensuremath{\operatorname L}\xspace}

   \providecommand{\Weaklyder}{\ensuremath{\operatorname W}\xspace}
   
   \providecommand{\dual}[1]{\ensuremath{{#1}'}}
   
   \providecommand{\dualspace}[2][]{\dual{\vecspace{#2}\ifx|#1|\else{_{#1}}}\fi}

   \providecommand{\cont}[1]{\ensuremath{\Continuous^{#1}}}

   \providecommand{\BV}[1]{\ensuremath{\operatorname{BV}}}

   \providecommand{\leb}[1]{\ensuremath{\Lebesgue_{#1}}}
   
   \providecommand{\lebnorm}[3][]{\ensuremath{\Norm{#2}_{\leb{#3}\ifx|#1|{}\else(#1)\fi}}}

   \providecommand{\sob}[2]{\ensuremath{{\smash\Weaklyder}^{#1}_{#2}}}
   
   \providecommand{\sobz}[2]{\ensuremath{{\overset{\smash{\scriptscriptstyle\circ}}\Weaklyder}{}^{#1}_{#2}}}
   
   \providecommand{\sobh}[1]{\ensuremath{\Hspace^{#1}}}
   
   \providecommand{\sobhz}[1]{\sobh{#1}_0}

   \providecommand{\Lip}[1][]{\ensuremath{\operatorname{Lip}}\ifx|#1|{}\else{\qp{#1}}\fi}

   \providecommand{\holder}[2]{\cont{#1,#2}}
   \providecommand{\poly}[1]{\ensuremath{\rP}^{#1}}

   \providecommand{\Symmatrices}[2][R]{\ensuremath{\operatorname{Sym}{(\csname#1\endcsname{#2})}}}\providecommand{\Symmats}[2][\reals]{\ensuremath{\operatorname{Sym}{({#1}^{#2})}}}
   \providecommand{\SAmatrices}[2][F]{\ensuremath{\operatorname{SA}{(\csname#1\endcsname{#2})}}}

  \providecommand{\crouzeixraviart}[1][1]{\operatorname{CR}\ifx|#1|{}\else{^{#1}}\fi}

  \providecommand{\linspace}[1]{\mathcal{\MakeUppercase{#1}}}
  \providecommand{\vecspace}{\linspace}

  \providecommand{\clinopss}[2]{\clinopss{\linspace{#1}}{\linspace{#2}}}
  \providecommand{\fepartition}[2][]{\mathscript{\MakeUppercase{#2}}\ifx|#1|{}\else_{#1}\fi}
  \providecommand{\fespace}[2][]{\mathbb{\uppercase{#2}}\ifx|#1|{}\else_{#1}\fi}
  \providecommand{\femesh}[2][]{\mathcal{\uppercase{#2}}\ifx|#1|{}\else_{#1}\fi}
  \providecommand{\vespace}{\fespace v}

  \providecommand{\fesh}[1]{\ensuremath{\fes{#1}h}}%

  \providecommand{\fe}[2][]{\ensuremath{\uppercase{#2}\ifx|#1|{}\else{_{#1}}\fi}}%
  \providecommand{\vecfe}[2][]{\ensuremath{\vec{\uppercase{#2}}\ifx|#1|{}\else{_{#1}}\fi}}%
  \providecommand{\matfe}[2][]{\ensuremath{\mat{\uppercase{#2}}\ifx|#1|{}\else{_{#1}}\fi}}%
  \providecommand{\hatmatfe}[2][]{\ensuremath{\hatmat{\uppercase{#2}}\ifx|#1|{}\else{_{#1}}\fi}}%
  
  \providecommand{\EOC}{\ensuremath{\operatorname{EOC}}\xspace}
  \providecommand{\tol}{\ensuremath{\operatorname{tol}}\xspace}

  \providecommand{\Foreach}{\text{ for each }}%

  \RequirePackage{amscd}

  \providecommand{\funk}[3]{\ensuremath{#1:#2\to#3}}
  
  \providecommand{\dfunkmapsto}[6][]{\ensuremath{
      \begin{array}{rrcl}
        {#2}: & {#4} &  \to   & {#6}
        \\
              & {#3} &\mapsto & {#5\text{\ #1}}
      \end{array}\quad}}

  \providecommand{\restriction}[2]{\left.#1\right|_{#2}}
  \renewcommand{\restriction}[2]{\left.#1\right|_{#2}}

  \providecommand{\evalat}[3][]{\qb{#2}_{\ifx|#1|{}\else#1=\fi#3}}
  \providecommand{\evaldiff}[4][]{\qb{#2}^{\ifx|#1|{}\else#1=\fi#3}_{\ifx|#1|{}\else#1=\fi#4}}

  \providecommand{\aka}[1]{(also known as {#1})\xspace}

  \providecommand{\codename}[1]{\url{#1}\xspace}
  \providecommand{\colorvarname}[2][]{\texttt{\ifx|#1|\else\color{#1}\fi#2}\xspace}

  \ifthenelse{\boolean{useutopia}}{%
    
  }{%
    
  }
  \providecommand{\indexen}[2][]{{\ifthenelse{\boolean{shownotes}}{\color a}{}#2\ifx|#1|\index{#2}\else\index{#1}\fi}}
  \providecommand{\indexemph}[1]{\emph{\indexen{#1}}}

  \providecommand{\eqncomment}[1]{\ensuremath{\qquad\qp{\text{\tiny{#1}}}}}

  \providecommand{\ListParameters}{}
  \renewcommand{\ListParameters}%
  {
  	 \setlength{\topsep}{0pt}
  	 \setlength{\leftmargin}{0pt}
           \setlength{\itemsep}{0pt}
  	 \setlength{\parsep}{0pt}
  	 \setlength{\parskip}{0pt}
           \setlength{\labelsep}{0pt}
  	 \setlength{\itemindent}{0pt}
  }
  {%
    \begin{list}%
      {}%
      {\ListParameters%
      
  }}%
  {\end{list}}
  \newcounter{tmpcounter}
  \newcounter{LetterListItem}
  \renewcommand{\theLetterListItem}{(\alph{LetterListItem})}
  {
  	\begin{list}%
  	{\theLetterListItem\ }%
  	{\usecounter{LetterListItem}
  	  \ListParameters
            \ifx|#1|{}\else\setcounter{LetterListItem}{#1}\fi
  	}
  }%
  {\end{list}}
  \newcounter{NumberListItem}
  \renewcommand{\theNumberListItem}{\arabic{NumberListItem}}
  {
  	\begin{list}%
  	{\theNumberListItem.\ }%
  	{\usecounter{NumberListItem}%
  	 \ListParameters
  	}
  }%
  {\end{list}}
  \newcounter{QuestionListItem}
  \renewcommand{\theQuestionListItem}{\textbf{Question \arabic{QuestionListItem}}}
  {
  	\begin{list}%
  	{\theQuestionListItem.\ }%
  	{\usecounter{QuestionListItem}%
  	 \ListParameters
  	}
  }%
  {\end{list}}
  \newcounter{RomanListItem}
  \renewcommand{\theRomanListItem}{(\roman{RomanListItem})}
  {
  	\begin{list}%
  	{\theRomanListItem\ }%
  	{\usecounter{RomanListItem}
  	 \ListParameters
  	}
  }%
  {\end{list}}
  \newcounter{StepsItem}
  {
  	\begin{list}%
  	{Step \theStepsItem.\ }%
  	{\usecounter{StepsItem}%
  	 \ListParameters
  	}
  }%
  {\end{list}}
  \newcounter{QAListItem}
  \renewcommand{\theQAListItem}{Q\arabic{QAListItem}:}
  {
  	\begin{list}%
  	{\theQAListItem}%
  	{\usecounter{QAListItem}
  	 \ListParameters
  	}
  }%
  {\end{list}}

  \providecommand{\grad}{}
  \renewcommand{\grad}[1][]{\nabla\ifx|#1|\else_{#1}\fi}

  \ifthenelse{\boolean{isthesis}}{%
    \setcounter{secnumdepth}{1}%
  }{
    \setcounter{secnumdepth}{2} %
  }
  \providecommand{\ListParameters}{}
  \renewcommand{\ListParameters}
  {
  	 \setlength{\topsep}{0em}
  	 \setlength{\leftmargin}{0em}
           \setlength{\itemsep}{0ex}
  	 \setlength{\parsep}{.5ex}
  	 \setlength{\itemindent}{\labelsep}
  	 \addtolength{\itemindent}{\labelwidth}
  }

    \providecommand{\ObsName}{Remark}%
    \providecommand{\RemName}{Remark}%
    \providecommand{\NotName}{Notation}%
    \providecommand{\BFNName}{Big~Fat~Note}%
    \providecommand{\DefName}{Definition}%
    \providecommand{\ExaName}{Example}%
    \providecommand{\TheName}{Theorem}%
    \providecommand{\LemName}{Lemma}%
    \providecommand{\ProName}{Proposition}%
    \providecommand{\CorName}{Corollary}%
    \providecommand{\PbmName}{Problem}%
    \providecommand{\HypName}{Hypothesis}%
    \providecommand{\AlgName}{Algorithm}%
    \providecommand{\ExeName}{Exercise}%
    \providecommand{\SolName}{Solution}%
    \providecommand{\ClaName}{Claim}%
    \providecommand{\EsyName}{Essay}%
    \providecommand{\Proofname}{Proof}%
    \providecommand{\Derivename}{Derivation}%
  
  \ifthenelse{\boolean{isthesis}}{%
    \providecommand{\Thecounter}{The}
  }{%
    \providecommand{\Thecounter}{subsection}
  }
  \makeatletter
  \newcommand{\oltikzgetxy}[3]{%
    \tikz@scan@one@point\pgfutil@firstofone#1\relax
    \edef#2{\the\pgf@x}%
    \edef#3{\the\pgf@y}%
  }
  \makeatother
  \providecommand{\pdfformat}[1]{
     \provideboolean{pdfoutput}
     \setboolean{pdfoutput}{#1}%
    \ifthenelse{\boolean{pdfoutput}}{
      \typeout{using pdf}
\usepackage{pdfsync}
      \providecommand{\graphext}{pdf}
      \renewcommand{\graphext}{pdf}
      \providecommand{\graphextex}{pdf_t}
      \renewcommand{\graphextex}{pdf_t}
    }{
      \typeout{using eps}
      \RequirePackage[dvips]{graphicx,xcolor}
      \providecommand{\graphext}{eps}
      \renewcommand{\graphext}{eps}
      \providecommand{\graphextex}{eps_t}
      \renewcommand{\graphextex}{eps_t}
    }
    \RequirePackage{epsfig}
    \RequirePackage{tikz}
    \RequirePackage{rotating}
    \RequirePackage{graphicx}
    \RequirePackage{xcolor}
    \provideboolean{darkcolortheme}
    \definecolor{SussexFlint}{rgb}{.00,.19,.21}
    \definecolor{SussexGrey}{rgb}{.51,.58,.49}
    \definecolor{SussexOrange}{rgb}{.94,.29,.00}
    \definecolor{SussexYellow}{rgb}{1.00,.73,.00}
    \definecolor{SussexRed}{rgb}{.94,.01,.49}
    \definecolor{SussexPurple}{rgb}{.48,.06,.44}
    \definecolor{SussexGreen}{rgb}{.00,.58,.46}
    \definecolor{OmarGreen}{rgb}{.00,.68,.36}
    \definecolor{SussexBlue}{rgb}{.00,.58,.65}
    \definecolor{OmarBlue}{rgb}{.00,.38,.65}
    \colorlet{a}{OmarBlue}%
    \colorlet{b}{SussexYellow}
    \colorlet{c}{SussexOrange}
    \colorlet{d}{OmarGreen}%
    \colorlet{e}{SussexRed}
    \colorlet{f}{SussexPurple}
    \colorlet{g}{white}%
    \colorlet{h}{SussexGrey}%
    \colorlet{i}{black}%
    \colorlet{j}{SussexFlint}
    \newcommand{\mausDarkColorTheme}{
      \colorlet{a}{SussexYellow!50!yellow}
      \colorlet{b}{SussexBlue}%
      \colorlet{c}{SussexOrange!50!yellow}
      \colorlet{d}{SussexGreen!50!green}
      \colorlet{e}{SussexRed!50!red}
      \colorlet{f}{SussexPurple!50!magenta}
      \colorlet{g}{black}%
      \colorlet{h}{SussexFlint!50!black}
      \colorlet{i}{white}%
      \colorlet{j}{SussexGrey}
    }
    \ifthenelse{\boolean{darkcolortheme}}{\mausDarkColorTheme}{}
  }
  \providecommand{\solution}{\textbf{\SolName.}\xspace}

   \newcounter{phantombox}[enumi]%
   \provideboolean{showphantoms}
   \renewcommand{\thephantombox}{\ifnum\value{phantombox}<10{0}\else{}\fi\arabic{phantombox}}
   \newcommand{\phantombox}[1]{\stepcounter{phantombox}%
     \ensuremath{\boxed{%
         {\ifthenelse{\boolean{showphantoms}}{#1}{\phantom{#1}}}%
         {\texttt{\tiny[\thephantombox]}
         }%
       }%
     }%
   }

   \provideboolean{hidesolution}
   \newcommand{\consolution}[2][]{
     \ifthenelse{\boolean{hidesolution}}{#1\setboolean{showphantoms}{false}}{%
       {\setboolean{showphantoms}{true}\color{i!50}\par \small {\solution}\ #2\par\ \\[5pt]}}
   }
   \provideboolean{showmarks}
   \providecommand{\showmarks}[1]{%
     \ifthenelse{%
       \boolean{showmarks}}{%
       \marginpar{%
         \tiny [$#1$ mark\ifthenelse{\equal{#1}1}{\phantom{s}}s]}%
     }{}}%

   \newcommand{\condibreak}{\ifthenelse{\boolean{hidesolution}}{\newpage}{}}

   \providecommand{\qeyword}[1]{\index{#1}\ifthenelse{\boolean{shownotes}}{\texttt{\color{e}[#1]}}{}}
   \providecommand{\sourcecite}[2][]{\index{#1}\ifthenelse{\boolean{shownotes}}{\texttt{\color{d}[source: \cite[#1]{#2}]}}{}}

   \RequirePackage{hyphenat}
\RequirePackage{xcolor,graphicx}
\RequirePackage{epsfig}
\RequirePackage{tikz}
\colorlet{a}{red!90!yellow}
\colorlet{b}{blue}
\colorlet{c}{green!50!blue}
\colorlet{d}{magenta}
\colorlet{e}{cyan}
\colorlet{f}{yellow!70!red}
\colorlet{g}{white}
\colorlet{h}{black!50}
\colorlet{i}{black}
\colorlet{j}{black!75} 

\newtheoremstyle{plain}%
  {}%
  {}%
  {\mdseries\slshape}%
  {\parindent}%
  {\bfseries}%
  {.}%
  {.5em}%
  {}%

\newtheoremstyle{note}%
  {}%
  {}%
  {}%
  {\parindent}%
  {\bfseries}%
  {.}%
  {.5em}%
  {}%

\newtheoremstyle{claim}%
  {}%
  {}%
  {\mdseries\slshape}%
  {}%
  {\bfseries}%
  {}%
  {.5em}%
  {}%

\newtheoremstyle{exercise}%
  {}%
  {}%
  {}%
  {}%
  {\bfseries}%
  {.}%
  {1em}%
  {}%

\newtheoremstyle{break}%
  {}%
  {}%
  {}%
  {}%
  {\bfseries}%
  {.}%
  {\newline}%
  {}%

\swapnumbers{
  \theoremstyle{plain}
  \ifthenelse
      {\boolean{isthesis}}
      {\newtheorem{The}{\TheName}[section]}%
      {
      }%
{
   \theoremstyle{plain}

   \renewcommand{\Thecounter}{subsection}

   \newtheorem*{The*}{\TheName}
   \newtheorem*{Lem*}{\LemName}
   \newtheorem*{Pro*}{\ProName}
   \newtheorem*{Cor*}{\CorName}
   \newtheorem*{Pbm*}{\PbmName}
   \newtheorem*{Hyp*}{\HypName}
   \newtheorem*{Exe*}{\ExeName}
   \newtheorem*{Txx*}{\ExeName} %
   \newtheorem*{Con*}{Conclusion}
   \newtheorem*{Sum*}{Summary}
 }
 {
   \theoremstyle{claim}

 }
 {
   \theoremstyle{note}

   \newtheorem*{Obs*}{\ObsName}

   \newtheorem*{Def*}{\DefName}
   \newtheorem*{Exa*}{\ExaName}
   \newtheorem*{Alg*}{\AlgName}
 }

 {
   \theoremstyle{break}
 }
}

\newenvironment{Lem}[1][]{\subsection{\LemName\xspace{\ifx&#1&{}\else{ (#1)}\fi}}\slshape}{\upshape}

\newenvironment{Def}[1][]{\subsection{\DefName\xspace{\ifx&#1&{}\else{ of #1}\fi}}}{}
\newenvironment{Obs}[1][]{\subsection{\ObsName\xspace{\ifx&#1&{}\else{ (#1)}\fi}}}{}

\newenvironment{Alg}[1][]{\subsection{\AlgName\xspace{\ifx&#1&{}\else{ (#1)}\fi}}}{}
\providecommand{\qed}{\vrule height 5pt depth 0pt width 3pt}
\providecommand{\qqed}{{\raggedright{\ \hfill\qed}}}

\newcounter{passo}

\newenvironment{Proof}[1][{}]%
{\par\noindent{\bf \Proofname\ #1}\setcounter{passo}{0}}%
{\qqed\par}
{\par\noindent{\bf \Derivename\ #1}\setcounter{passo}{0}}%
{\qqed\par}
\newenvironment{Proof*}[1][{}]%
{\subsection{\Proofname\ #1}\setcounter{passo}{0}}
{\qqed\par}

\usepackage{natbib}
\usepackage{algpseudocode}
\providecommand{\fenics}{FEniCS\xspace}
\newcommand{\mypart}[1][]{\fepartition T\ifx|#1|\else_{#1}\fi}
\newcommand{\Wh}[1][]{\W_{\mypart[#1]}}
\providecommand{\det}{\operatorname{det}}
\providecommand{\itermax}{\operatorname{itermax}}
\providecommand{\tol}{\operatorname{tol}}

\providecommand{\Cof}{\operatorname{Cof}}
\providecommand{\xrp}{\vec x,r,\vec p}
\providecommand{\xp}{\vec x,\vec p}
\renewcommand{\fesh}{\fespace}
\providecommand{\feop}[2][]{\mathsf{\MakeUppercase{#2}}\ifx|#1|{}\else_{#1}\fi}
\providecommand{\vecfeop}[2][]{\vec{\mathsf{\MakeUppercase{#2}}}\ifx|#1|{}\else_{#1}\fi}
\providecommand{\matfeop}[2][]{\mat{\mathsf{\MakeUppercase{#2}}}\ifx|#1|{}\else_{#1}\fi}
\providecommand{\tildematfeop}[2][]{\tilde{\matfeop{#2}}\ifx|#1|{}\else_{#1}\fi}
\providecommand{\lproj}[1][h]{\feop P}
\usepackage{caption}
\usepackage{subcaption}
\usepackage{rotating}
\usepackage{nicefrac}
\setboolean{shownotes}{false}%
\setboolean{showchanges}{false}%
\setboolean{usemathrsfs}{false}%
\makeindex
\author{Ellya Kawecki}
\thanks{EK acknowledges support of the Engineering and Physical Sciences Research Council [EP/L015811/1]}
\author{Omar Lakkis}
\author{Tristan Pryer}
\title[Finite elements for Monge--Ampère with optimal transport] {
  A finite element method for the Monge--Ampère equation with
  optimal transport boundary conditions
}
\date{\today}
\begin{document}
\maketitle
\begin{abstract}
  We address the numerical solution via Galerkin type methods of the
  Monge--Ampère equation with optimal transport boundary conditions, arising in
  optimal mass transport, geometric optics and mesh/grid movement
  techniques.  This fully nonlinear elliptic problem admits a
  linearisation via a Newton--Raphson iteration, which leads to a sequence of
   elliptic equations in
  nondivergence form, with oblique derivative boundary conditions.  We discretise these by employing the
  nonvariational finite element method, which leads to empirically
  observed optimal convergence rates, provided recovery techinques are
  used to approximate the gradient and the Hessian of the unknown
  functions.  We provide extensive numerical testing to illustrate the
  strengths of our approach and the potential applications in optics
  and mesh movement.
 \end{abstract}
\section{Introduction}
\subsection{The Monge--Ampère problem}
\label{sec:monge-ampere-problem}
Given $d\geq1$, two convex domains (i.e., open and bounded subsets) of
$\R d$, $\W $ and $\Y $, and two uniformly positive functions
$\funk\rho\W\reals$ and $\funk\sigma\Y\reals$ (thought as mass densities)
with equal total mass ($\smallint_\W {\rho}=\smallint_\Y {\sigma}$),
the \indexemph{Monge--Ampère problem for optimal transport} \emph{(\indexen{MAOT})}
consists of finding a function
$\funk u\W\reals$ satisfying the following \indexemph{domain transport condition}
\begin{equation}
  \label{eq:monge-ampere-domain-transport-condition}
  \grad u(\W )=\Y
  ,
\end{equation}
and the \indexemph{partial differential equation} (\indexemph{PDE})
\begin{equation}\label{eq:monge-ampere-problem}
  \det\D^2 u(\vec x)=\frac{\rho(\vec x)}{\sigma(\grad u(\vec x))}
  \text{ for }\vec x\in\W 
  ,
\end{equation}
where $\grad u(\vec x)$ and $\D^2 u(\vec x)$ respectively denote the
gradient vector and Hessian tensor (or matrix) of $u$ at $\vec x$.

The PDE (\ref{eq:monge-ampere-problem}) is commonly known as the
\indexemph{Monge--Ampère equation} and the domain transport condition in
(\ref{eq:monge-ampere-domain-transport-condition}) is called the
\indexemph{transport boundary condition} \aka{\indexemph{second boundary condition}}
because with (\ref{eq:monge-ampere-problem}) it defines a boundary value problem
associated to the optimal transport problem of finding a map
$\funk{\vec t}\W {\R d}$, known as the \indexemph{transport map} or
\indexemph{transport field}, such that
\begin{equation}
  \label{eq:Monge-Ampere-Jacobi-problem}
  \vec t(\W )=\Y 
  \tand
  \sigma(\vec t(\vec x))\det \D\vec t(\vec x)=\rho(\vec x).
\end{equation}
Here the equation can be interpreted as the change of variables \(\vec
y=\vec t(\vec x)\), which clarifies why we required equality of mass
for $\rho$ and $\sigma$ as
\begin{gather}
  \d\vec y=\det\D\vec t(\vec x)\d\vec x
  \intertext{and thus }
  \int_\W  \rho(\vec x)\d\vec x
  =
  \int_\Y  \sigma(\vec t(\vec x))\det \D\vec t(\vec x)\d\vec x
  =
  \int_\Y  \sigma(\vec y)\d\vec y
  .
\end{gather}
At first sight it may appear
that~(\ref{eq:monge-ampere-domain-transport-condition}) is not a boundary
condition, since values of the gradient are prescribed in the
\emph{interior} of $\W$ instead of its boundary, $\boundary\W$.  But
on closer inspection, and as shown in
\citet{Urbas:97:article:On-the-second}, this condition, thanks to the
convexity of the domains $\W$ and $\Y$, is in fact equivalent to the
boundary-only condition
\begin{equation}
  \vec t(\boundary\W)=\boundary\Y.
\end{equation}

Thanks to the polar factorization of transport maps discovered by
\citet{Brenier:91:article:Polar}, a field $\vec t$, with simply
connected uniformly convex domains $\W$ and $\Y$, satisfying
(\ref{eq:Monge-Ampere-Jacobi-problem}), exists if and only if there
exists a uniformly convex function $\funk u\W\reals$, such that 
$\vec t=\grad u$, satisfying~(\ref{eq:monge-ampere-problem})
and~(\ref{eq:monge-ampere-domain-transport-condition}); and the latter can be
replaced by
\begin{equation}\label{eq:monge-ampere-realboundaryproblem}
  \grad u(\boundary\W )=\boundary\Y .
\end{equation}
In
\citet{PrinsTen-Thije-BoonkkampRoosmalenIjzermanTukker:14:article:A-Monge-Ampere-solver}
the equivalence is proved under the weaker assumption that $\W $ and
$\Y $ are simply connected domains, the existence and uniqueness
results proven by \citet{Urbas:97:article:On-the-second} require
that both $\W $ and $\Y $ are uniformly convex. As such this will be
our general assumption unless stated otherwise.  An overview can
also be found in
\citet{PrinsTen-Thije-BoonkkampRoosmalenIjzermanTukker:14:article:A-Monge-Ampere-solver}
and more in \citet[Ch.4]{Villani:03:book:Topics}.

For $1=d$ the Monge--Ampère equation reduces to the textbook
(linear) Poisson equation and the transport boundary condition to
a nonlinear Neumann boundary condition. For $2\leq d$, equation
(\ref{eq:monge-ampere-problem}) is a fully nonlinear elliptic equation,
while the boundary condition (\ref{eq:monge-ampere-domain-transport-condition})
equation can be interpreted as a nonlinear condition on the gradient
(or a Hamilton--Jacobi equation) of the function $\grad u$.

MAOT arises in many areas of
mathematics, such as differential geometry, meteorology, and the
design of free form reflectors.  One particular meteorological
application of the MAOT problem is the incorporation of moving meshes
in the solution of meteorological partial differential
equations. \citet{BuddCullenWalsh:13:article:MongeAmpere} successfully coupled a parabolic MAOT
method for the construction of a moving mesh in two-dimensions to a
pressure correction method. In this case the MAOT problem serves to
generate a moving mesh on which the PDE is solved numerically
\citep{BuddRussellWalsh:15:article:The-geometry}.

The linearisation of MA type equations typically results in a sequence
of nondivergence form elliptic equations. Such problems do not, in
general, possess a weak formulation, and as such, standard conforming
finite element methods must be either restricted
\citep{NochettoZhang:18:article:Discrete}, modified into mixed forms
\citep{LakkisPryer:11:article:A-finite,Gallistl:17:article:Variational},
nonconforming (\emph{discontinuous}
Galerkin)~\citep{SmearsSuli:13:article:Discontinuous,Kawecki:17:techreport:A-DGFEM-oblique},
or obtained in the limit of fourth-order
perturbations~\citep{FengNeilan:14:article:Finite}.
\citet{LakkisPryer:13:article:A-finite,LakkisPryer:11:article:A-finite}
proposed a \emph{continuous} Galerkin finite element method, called
the nonvariational finite element method (NVFEM), which approximates
solutions of nondivergence form elliptic problems, with Dirichlet
boundary conditions. Other notable \emph{discontinuous} Galerkin
finite element methods were derived
by~\citet{SmearsSuli:13:article:Discontinuous} in the context of
convex polytopal domains, as well
as~\cite{Kawecki:17:techreport:A-DGFEM} in the context of curved
domains with piecewise nonnegative curvature.  It is worth noting an
alternative approach using semilagrangian methods on Galerkin-type
(and therefore not necessarily structured) meshes by
\citet{FengJensen:17:article:Convergent} which has the potential to be
exported to optimal transport conditions.  For more in depth
information about the state of the art on numerical methods for
Monge--Ampère type PDEs and related boundary value problems, we refer
to the review of \citet{NeilanSalgadoZhang:17:article:Numerical}.
\subsection{Literature and context}
\label{landc:subsec}
Assuming $\W$ and $\Y$ are uniformly convex $\holder21$ domains,
$\rho,\sigma\in\holder11(\closure\W)$,
\citet{Urbas:97:article:On-the-second} proves the existence of a
convex function
$u\in\holder3\alpha(\W)\meet\holder2\alpha(\closure\W)$, for all
$\alpha\in\opinter01$ satisfying (\ref{eq:theprob}), as well as its
uniqueness up to an additive constant.
\citet{Urbas:97:article:On-the-second}'s
idea is to represent the target domain $\Y $ as the superlevel set of a concave
defining function $\funk b{\R d}\reals$, i.e.,
\begin{equation}
  \Y =\{\vec p\in\R d:b(\vec p)>0\}.
\end{equation}
It can then be seen that $\boundary\Y =\{\vec p\in\R{d}:b(\vec
p)=0\}$.  The problem can then be recast in the following
\indexemph{nonlinear second boundary value elliptic problem}
\begin{gather}
  \label{eq:theprob}
  \begin{split}
    \det\D^2u(\vec x)
    &
    =\frac{\rho(\vec x)}{\sigma(\nabla u(\vec x))}\quad\mbox{for}\quad\vec x\in\W,
    \\
    b(\nabla u(\vec x))
    &
    =0\quad\mbox{for}\quad\vec x\in\boundary\W.
  \end{split}
\end{gather}

Using this formulation,
\citet{BenamouFroeseOberman:14:article:Numerical} provided a numerical
method based on the \indexemph{wide-stencil finite difference}
approach and a treatment of the boundary via a Hamilton--Jacobi
approximation.  The scheme they provide is consistent and monotone in
the sense of \citet{BarlesSouganidis:91:article:Convergence} (and thus
convergent) and numerical experiments show that it cannot be more than
first order, which is to be expected for such monotone
schemes. \citet{BenamouFroeseOberman:14:article:Numerical} state that
the accuracy can be somewhat restored by making the schemes ``almost
monotone'' (sic) referring to \citet{Abgrall:09:article:Construction}
without giving much details; it is also unclear how monotonicity is
ensured when second boundary conditions are prescribed as opposed to Dirichlet
boundary conditions. To numerically encode the boundary
condition~(\ref{eq:monge-ampere-domain-transport-condition}),
\citet{BenamouFroeseOberman:14:article:Numerical} used a similar
representation to what is seen in
\citet{Urbas:97:article:On-the-second} with a \emph{convex} (instead
of concave) defining function $b:\R{d}\to\reals$, so that $\Y
=\setofsuch{\vec p\in\R d}{b(\vec p)<0}$.  In particular they pick for
$b$ (which is not unique) the signed distance function of the target
boundary $\Y$, that is
\begin{equation}
  \label{bdef}
  b(\vec p)=
  \qp{1-2\iverson{\vec p\in\Y}}
    \dist(\vec p,\boundary\Y )
    ,
\end{equation}
where, for a proposition $P$, the \indexemph{Iverson--Knuth bracket} is
\index{$\iverson\cdot$}$\iverson P:=1$ if $P$ is true, $0$ if $P$ is false.
\subsection{Our main results}
In this article, we propose a finite element method that works for
$\poly k$ (polynomial of degree $k\in\naturals$) elements for any $k$
and on unstructured meshes which leads to convergence with high order,
in many cases optimal, and opens the way to adaptive mesh refinement
strategies for problems with singular, e.g., viscosity, solutions in
the spirit of \citet{Pryer:10:phdthesis:Recovery},
\citet{LakkisPryer:15:inproceedings:An-adaptive} and
\citet{Gallistl:17:article:Variational}.
Our method consists in
\begin{enumerate}
\item
  introducing a Newton--Raphson's method at the continuum
  (exact) stage that iteratively approximates the solution of
  \cite{LakkisPryer:15:inproceedings:An-adaptive} (this results
  in a sequence of \indexemph{oblique derivative} boundary value problems
  for elliptic equations in \indexemph{nondivergence form});
\item
  applying a nonvariational finite element method (NVFEM) similar to the
  one from \cite{LakkisPryer:13:article:A-finite} but including the
  oblique boundary condition for which gradient recovery techniques
  prove crucial in order to obtain (empirically observed) optimal
  convergence rates (incidentally, the gradient recovery is also
  useful to provide a finer approximation of the argument of the
  target density, $\sigma$).
\end{enumerate}
We apply a global gradient recovery scheme, similar to the one
described in \citet{ZhangNaga:05:article:A-new-finite} and
\citet{ZienkiewiczZhu:87:article:A-simple}, in order to achieve
convergence of the algorithm for $\poly1$ elements. Without the use of
gradient recovery the scheme is seen to only converge for $\poly{k}$
elements, where $k\ge 2$.  While finalising this paper, a closely
related one by \citet{Gallistl:17:techreport:Numerical} was brought to
our attention; therein the author tackles the oblique derivative
problem with mixed-method techniques.

The rest of this paper is organised as follows:
In \S\ref{NandF:sec} we provide the notation needed, and define the
finite element spaces we use in our numerical method. In \S\ref{nvfem}
we introduce the nonvariational finite element method for problems
that arise in the linearisation of~(\ref{eq:theprob}). In
Section~\ref{oblique:subsec}, further details on how we adapt the
method found in~\citet{LakkisPryer:13:article:A-finite,LakkisPryer:11:article:A-finite}
to the context of oblique boundary value problems will be given.  In
\S\ref{sec:nonlinear-solver} we discuss the conditional
ellipticity of the nonlinear operator associated
to~(\ref{eq:monge-ampere-problem}), and appropriate linearisation
schemes.  In \S\ref{sec:The scheme} we define the two numerical
methods that are the main focus of this paper, the second method is
distinguished from the first by the inclusion of a gradient recovery
operator in the finite element scheme.  We report on our numerical experiments
in \S\ref{results}, by looking first at cases where the true solution is
known, so that we can observe the rates of convergence of the
numerical methods, followed by experiments where the true solution is
unknown, demonstrating the robustness of the two methods.  Finally, in
\S\ref{sec:conclusion} we give concluding remarks on what has
been accomplished in this paper, as well as plans for future research.
\section{Notation and functional set up}
\label{NandF:sec}%
\subsection{Vector, matrix and function spaces}
The usual real $d$-dimensional Euclidean space is denoted $\R d$ with
$\norm\cdot$ denoting the (Euclidean) norm.  We write $\realmats {m}{n}$
for the space of all real coefficient matrices with $m$ rows and $n$
columns; this is identified with the space of linear transformations
from $\R n$ into $\R m$. The \indexemph{Frobenius product} \aka{\indexemph{double dot product}}
of two matrices, say 
$\mat{A}=\matijnm{a}{i}{j}{m}{n}$ and $\mat{B}=\matijnm{b}{i}{j}{m}{n}$, 
in \(\realmats {m}{n}\) is defined as\index{$:$}
\begin{equation}
  \label{eq:def:Frobenius-inner-product}
  \mat A\frobinner\mat B
  :=
  \trace\qp{\transposemat A\mat B}
\end{equation}
where $\transposemat A$ is matrix $\mat A$'s \indexen{transpose matrix}
and $\trace\mat M$ is matrix $\mat M$'s \indexen{trace}.
Immediate properties of the Frobenius product are
\begin{equation}
  \mat A\frobinner\mat B
  =
  \sumifromto{i,j}1{m,n}\matentry aij\matentry b i j
  \tand
  \mat A\frobinner\mat B
  =
  \trace\qp{\mat A\transposemat B}
  .
\end{equation}
The Frobenius product turns the space of linear operators $\realmats
{m}{n}$ into a Hilbert space, which, for $n=1$ (or $m=1$), trivially
coincides with the usual Euclidean space of vectors (or
covectors).  The set of all \indexen{symmetric} operators (matrices)
\indexen{$\Symmatrices d$} is a
linear subspace of $\realmats dd$.  The set of all \indexen{symmetric and
positive definite} operators (matrices) on $\R d$, \indexen{$\spdmats[R] d$},
is a subset of $\Symmatrices d$.

Wherever we use measure and integration we intend Lebesgue's, if the
integration domain has non-zero Lebesgue measure (with elementary
measure $\d x$) or the surface, line, point (Hausdorff) measure.  with
and we indicate the surface measure \aka{$d-1$-dimensional Hausdorff
  measure} with \indexen{$\measure S$}.  We also omit the integration
elements, $\d\vec x$ or $\ds(\vec x)$, wherever the integration
variable $\vec x$ is silent or the meaning of the measure obvious from the
integration domain.

Let $K$ be an open or closed (Lebesgue or Hausdorff) measurable subset
of $\R d$. We consider the well-known spaces of $p$-summable
functions, for any real number $p\geq1$, \index{$\leb p$!$(p\geq1)$}
\begin{equation}
  \label{Lpspace}
  \leb p(K):=\ensemble{\funk v K\reals}{\int_K|v|^p<\infty},
\end{equation}
also defined when $p=\infty$ as\index{$\leb\infty$}
\begin{equation}
  \label{Linfspace} 
  \leb\infty(K)
  :=\ensemble{\funk v K\reals}{
    \exists M\in\R+:
    \abs{v(\vec x)}
    \leq M
    \Aa\vec x\in K
  }
  .
\end{equation}
We equip the spaces
$\leb p(K)$, $1\le p<\infty$ and $\leb \infty(K)$ with the following
norms
\begin{equation}
  \Norm{v}_{\leb p(K)}
  :=
  \powqp{\fracl1p}{\int_K\powabs p v}
  ,
\end{equation}
\begin{equation}
  \Norm v_{\leb\infty(K)}
  :=
  \inf\ensemble{M\in\R+}{|v(\vec x)|\le M\Aa\vec x\in K},
\end{equation} 
respectively. The space $\leb p(K)$ is a Banach space for any
$p\in\clinter1\infty$ \citep{LiebLoss:01:book:Analysis}.
  
In the special cases $p=1$ and $p=2$, we equip $\leb p(K)$,
respectively, with the continuous linear functional, and inner
product, respectively denoted\index{$\qa\cdot$}
\index{$\ltwop\cdot\cdot$}
\begin{equation}
  \qa{w}_K
  :=
  \int_K w
  \tand
  \ltwop u v_K
  :=
  \qa{u v}_K
  =
  \int_K u v
  \text{ for }
  w\in\leb1(\W),
  u,v\in\leb2(\W).
\end{equation}
The same notations are used also for tensor (including vector) valued
functions with the result returning a real valued tensor (vector) or
a scalar depending on the context.  The space $\leb{2}(K)$ is a Hilbert
space when equipped with $\ltwop\cdot\cdot_K$.  More generally,
whenever $\cV$ is a topological vector space and $\dual\cV$ its dual,
we indicate the duality pairing with
\begin{equation}
  \duality lv_\cV
  \text{ for }
  l\in\dual\cV,
  v\in\cV
  .
\end{equation}
For the sake of presentation, we often drop the subindex $\cV$, unless there is any chance of ambiguity.

We denote by $\D u$ the (possibly only distributional) derivative of a function
$\funk u K\R d$, and we define the gradient of $u$, $\grad u$, to be
the derivative's transpose, i.e.,\index{$\leb p$}
\begin{equation}
  \grad u=\transposeof{\D u}.
\end{equation}
For any $m\in\NO$, the $m$-th (possibly only distributional) derivative of $u$
is recursively defined by $$\D^m u:=\D\D^{m-1}u\,\,\mbox{and}\,\,\D^0u:=u.$$  We
denote by $\D^2u$ interchangeably the second derivative and the
Hessian of $u$, i.e., the $\Symmatrices d$ matrix of second order partial
derivatives of $u$; we prefer this abuse of notation to the more
consistent yet cumbersome notation for the Hessian as derivative of the
gradient, $\grad\D u$.

For $m\in\NO$ we introduce the following
Sobolev spaces~\citep{Evans:10:book:Partial}:\index{$\sob mp$}
\begin{equation}
  \begin{gathered}
    \sob mp(K):
    =
    \ensemble{v\in \leb p(K)}{
      \D^{\vec\alpha}v\in \leb 
      p(K)\Foreach\vec\alpha:|\vec\alpha|\le m},
    \text{ for }1\le p\le\infty,
    \\
    \text{ and the shorter form }
    \sobh m(K):=\sob m2(K),
  \end{gathered}
\end{equation}
using the multi-index notation
$\vec\alpha=\seqidotsfromto\alpha 1d\in\NO^d$, with
$|\vec\alpha|:=\sum_{i=1}^d\alpha_i$, and the partial derivatives,
$\D^\alpha=\pd 1^{\alpha_1}\dotsm\pd d^{\alpha_d}$, are understood in
the weak sense.  The spaces $\sob mp(K)$ are Banach with the following norms:
\begin{gather}
  \Norm{v}_{\sob mp(K)}
  :
  =\powqpBig{1/p}{
    \sum_{|\vec\alpha|\le m}\Norm{\D^{\vec\alpha} v}_{\leb p(K)}^p
  },
  \text{ if }
  1\le p\le\infty,
  \\
  \Norm{v}_{\sob m\infty(K)}:
  =\max_{|\vec\alpha|\le m}\Norm{\D^{\vec\alpha}v}_{\leb\infty(K)},
\end{gather}
and seminorms:
\begin{equation}
  \Norm v_{\sob mp(K)}:=\Norm{\D^mv}_{\leb p(K)}.
\end{equation}
The space $\sobh m(K)$ is Hilbert with the inner product
\begin{equation}
  \ltwop uv_{H^m(K)}
  :=
  \sum_{|\alpha|\le m}
  \int_KD^\alpha u~D^\alpha v.
\end{equation}
We define the space \(\sobz{m}p(K)
:=\setofsuch{v\in\sob{m}p(K)}{\restriction{v}{\boundary{K}}=0}\) where
a function's restriction to the boundary is understood as its \indexen{trace}
(not to be confused with a matrix's trace)
\citep{Evans:10:book:Partial}. We will occasionally use the fractional
order boundary Sobolev space
\indexen{$\sobh{\fracl12}(\boundary{K})$}, for an open $K$ of class
$\holder01$, to be thought of as the image of $\sobh1(K)$ under the
trace operator with the following norm:
\begin{equation}\label{H1/2normdef}
  \Norm{v}_{\sobh{\fracl12}(\boundary\W)}
  :=
  \inf
  \setofsuch{
    \Norm{w}_{\sobh1(\W)}
  }{
    w\in\sobh1(\W)
    \tand
    \restriction w{\boundary\W}=v
  }
  .
\end{equation}
\subsection{Finite element spaces}
The finite element spaces we consider will always be given with
respect to the domain $\W$ or a given subset of $\W$.

Consider \indexen{$\mypart$}
to be a \indexemph{fitted shape-regular triangulation} of $\W$,
namely $\mypart $ is a nonempty finite family of sets such that:
\begin{enumerate}[(i)\ ]
\item
  $K\in\mypart$ implies that $K$ is an open simplex (point for
  $d=0$, segment for $d=1$, triangle for $d=2$, tetrahedron for
  $d=3$);
\item 
  for any $K,J\in\mypart $ we have that $\closure{K}\meet\closure{J}$ 
  is a full closed subsimplex of $\closure K$ or $\closure J$.
\end{enumerate}
The triangulation's union may not coincide with $\W$, so we introduce the
\indexemph{approximate domain}\index{$\Wh$}
\begin{equation}
  \Wh:=\operatorname{int}\left(\cup_{K\in\mypart}\overline{K}\right)
\end{equation}
and note that (thanks to the convexity of $\W$) ${\Wh}\subseteq\W$.
We will also use the triangulations \indexemph{meshsize} (scalar)\index{$h$ (meshsize)}
\begin{equation}
  h:=\max_{K\in\mypart}\diam K\text{ with }
  \diam X:=\sup_{{\vec x},{\vec y}\in X}\norm{\vec x-\vec y}.
\end{equation}
All of our work can be replicated on more general partitions, involving
not only simplices but also other types of polytopes, but we do not
treat those in this paper to avoid distractions from our main goal.
Another, very useful generalisation would be the use of isoparametric
elements to approximate the boundary at an order higher than $2$, which is
what we presently do with straight elements.

We will use the following notation, valid for a generic vector
space $\cX(D;E)$ of functions with domain and range $D,E\subseteq\R d$, 
denoting by $O$ the union of all open simplices of $\mypart$ (note
that if $\mypart$ is not a singleton $O$ is a proper subset of $\W$)
\begin{equation}
  \cX(\mypart)
  :=
  \ensemble{\funk vOE}{\restriction vK\in\cX(K)\Foreach K\in\mypart}
\end{equation}
We say that elements of $\cX(\mypart)$ are piecewise (or $\mypart$-wise)
in $\cX$.
Let $\poly k$ denote the space of polynomials in $d$ variables of
degree less than or equal to $k\in\naturals$, and $\poly k(K)$ the
restriction of such functions to $K$; this allows us to define
the finite element spaces:
\begin{equation}\label{def:FEspaceV}
  \fesh v:=
  \poly k(\mypart )
  \meet
  \cont0(\overline{\Wh})
  =
  \ensemble{v\in\cont0(\overline{\Wh})}{\restriction vK\in\poly k(K)\Foreach K\in\mypart },
\end{equation}
as well as
\begin{equation}\label{def:FEspaceUW}
  \fesh g:=\powqp d{\fesh v},
  \tand
  \fesh h:=\setofsuch{
    \matfe W\in\powqp{d\times d}{\fesh v}
  }{
    \matfe W(\vec x)\in\Symmatrices d\Foreach\vec x\in\Wh}.
\end{equation}
The maximal polynomial degree $k\ge 1$ is fixed with respect to the
mesh elements, we denote by \index{$N$}$N:=\dim\fesh v$, the number of
the finite-element space's \indexemph{degrees of freedom} (\emph{\indexen{DOF}s}),
and an (ordered) nodal
basis \indexen{$\vecidotsfromto{\fe\varPhi}1N$} of $\fesh v$.
\section{The nonvariational finite element method}
\label{nvfem}
We now adapt the nonvariational finite element method (NVFEM) proposed
in~\citet{LakkisPryer:13:article:A-finite,LakkisPryer:11:article:A-finite}
to build a finite element approximation to $u$ satisfying
(\ref{eq:monge-ampere-problem}) and
(\ref{eq:monge-ampere-realboundaryproblem}).
\subsection{Linear nonvariational oblique derivative problem}
\label{oblique:subsec}
Let $\W$ be a convex ${\holder2{1}}$ domain, denote its outer normal, a unit vector-valued
function defined on \Aa[\measure s] of $\boundary\W$, by
$\normalto\W$. Let $\alpha\in(0,1)$ and $\mat A$ a symmetric uniformly
positive definite matrix-valued function in $\holder0\alpha(\W
;\realmats dd)$, i.e., there exists a constant $\mu>0$ such that
\begin{equation}\label{oblique:bvp}
  \vec\xi\transposed
  \mat A(\vec x)
  \vec\xi
  \geq
  \mu|\vec\xi|^2
  \Foreach 
  \vec\xi\in\R d,
  \vec x\in\W, 
\end{equation}
vector valued functions $\vec\beta\in\holder1\alpha(\boundary\W;\R d)$
such that for a constant $\beta_\flat>0$,
$\vec\beta\inner\normalto\W\geq\beta_\flat$ a.e.,
$\vec b\in\holder 0\alpha(\W;\R d)$, $c\in\holder0\alpha(\W)$, $c\le0$,
$r\in\holder0\alpha(\W)$, and $s\in\holder1\alpha(\boundary\W)$ find
$\funk u\W\reals$ that satisfies
\begin{equation}\label{genNDVobl}
  \begin{aligned}
    \mat A(\vec x):\D^2u(\vec x)
    +
    \vec b(\vec x)\inner\grad u(\vec x)
    +
    c(\vec x)u(\vec x)
    &
    =
    r(\vec x)
    \text{ for }
    \vec x\in\W
    ,
    \\
    \vec\beta(\vec x)\cdot 
    \grad u(\vec x)
    &=
    s(\vec x)
    \text{ for }
    \vec x\in\boundary\W.
  \end{aligned}
\end{equation}
The problem given above is an oblique derivative problem, which is
well posed in view of~\citet[Th.6.31,~e.g.]{GilbargTrudinger:01:book:Elliptic} when
$\W\in\holder2\alpha(\W)$ and
\citet{Lieberman:01:article:Pointwise,Lieberman:87:article:Local} for
Lipschitz domains, which comprise the herein needed convex ${\holder2{1}}$ domains.
\begin{Def}[generalised Hessian]\label{genhessdef}
To define the notion of the finite element Hessian, we must first
introduce the concept of the \emph{generalised Hessian}.  Looking
first at a smooth function, say
$v\in\cont2(\W)\meet\cont1(\closure\W)$, an application of integration
by parts shows us that the Hessian of $v$, $\D^2v$, satisfies
(the system of $d\times d$ equations)
\begin{equation}
  \label{eq:anappof}
  \ltwop{ \D^2v}\varphi
  =
  -
  \inta{{\grad v}{\D\varphi}}
  +
  \inton{\grad v\Transpose{\normalto\W}\varphi}{\boundary\W}
  \Foreach
  \varphi
  \in 
  \sobh 1(\W ),
\end{equation}
where $\normalto\W$ is the unit outward normal to $\W $.  We
generalise this to a given function $v\in \sobh 1(\W )$ with
$\restriction{\grad
  v\Transpose{\normalto\W}}{\boundary\W}\in\psetmats{\dual{\sobh{1/2}(\boundary\W)}}dd$ by
defining the \indexemph{generalised Hessian} of $v$, $\D^2v$, maps to
$\realmats dd$ rather than $\reals$ as an element in $\powqp{d\times
    d}{\dual{\sobh 1(\W)}}$ via
\begin{equation}
  \label{eq:def:generalised-Hessian}
  \duality{\D^2v}{\varphi}
  :=
  -
  \inta{
    \grad v
    \D\varphi
  }
  +
  \duality{
    \grad v
          {\normalto\W}\transposed
  }{
    \varphi
  }_{(\sobh{1/2}(\boundary\W))'\times\sobh{1/2}(\boundary\W)}
  \Foreach
  \varphi\in 
  \sobh 1(\W).
\end{equation}

Due to the duality pairing on the right-hand side of
(\ref{eq:def:generalised-Hessian}) our definition of generalised
Hessian is a $\setmats{\sobh1(\W)}dd$-continuous linear extension of the distributional Hessian 
to include test functions whose support needs not be compact in $\W$. Note
however that \emph{$\D^2v$ is not a distribution}.  Nevertheless, it is
a continuous linear functional, which legitimises our use of the duality
brackets $\duality\cdot\cdot$ to manipulate it.
\end{Def}
  \begin{Lem}[generalised Hessian linear functional]
    Assume that
    \begin{equation}
      v\in\sobh1(\W)\tand
      \grad{v}\normalto\W\transposed
      \in
      \psetmats{\dual{\sobh{1/2}(\boundary\W)}}dd
    \end{equation}
    then the right-hand side of
    (\ref{eq:def:generalised-Hessian})
    is a well-defined linear functional 
    \begin{equation}
      \D^2v\in\psetmats{\dual{\sobh1(\W)}}dd
      .
    \end{equation}
  \end{Lem}
  \begin{Proof}
    First
    \label{lem:generalised-Hessian-linear-functional}
    define the linear map
    $\funk{\mat W_v}{\sobh1(\W)}{\realmats dd}$
    \begin{equation}
      \duality{\mat W_v}{\varphi}
      :=
      -
      \inta{{\grad v}{\D\varphi}}
      +
      \duality{\grad v\Transpose{\normalto\W}}{\varphi}_{\sobh{1/2}(\boundary\W)}
      \text{ for }\varphi\in\sobh1(\W).
    \end{equation}
    Looking at each component of the resulting matrix, we see that
    for $i,j\integerbetween1d$,
    \begin{equation}
      \begin{split}
        \sup_{\Norm{\phi}_{\sobh1(\W)}=1}
        \duality{\getentryij{\mat W_v}ij}{\varphi}
        &
        \leq
        \sup_{\Norm{\varphi}_{\sobh1(\W)}=1}\ltwop{\pd i v}{\pd j\varphi}
        +
        \sup_{\Norm{\varphi}_{\sobh1(\W)}=1}\duality{\pd i v\getrow{\normalto\W}j}{\varphi}
        \\
        &
        \leq
        \sup_{\Norm{\varphi}_{\sobh1(\W)}=1}
        \qp{
          \Norm{v}_{\sobh1(\W)}\Norm{\varphi}_{\sobh1(\W)}
        }
        +
        \Norm{\pd iv\getrow{\normalto\W}j}_{\dual{\sobh{\fracl12}(\boundary\W)}}
        \\
        &\le 
        \Norm{v}_{\sobh1(\W)}
        +\Norm{
          \grad v{\normalto\W}\transposed
        }_{
          \psetmats{\dual{\sobh{\fracl12}(\boundary\W)}}dd}
        .
      \end{split}
    \end{equation}
    It then follows that \({
      \mat W_v\in\powqp{d\times{d}}{\dual{\sobh1(\W)}}
    }\), and thus the right hand side
    of~(\ref{eq:def:generalised-Hessian}) is well defined and the
    linear functional $\D^2v$ thus defined is in
    $\psetmats{\dual{\sobh1(\W)}}dd$.
  \end{Proof}
\begin{Def}[finite element Hessian]
  \label{femhess}
  From \eqref{eq:def:generalised-Hessian}, and in view of the Riesz
  representation theorem, for $v\in\fesh v$ we define its ($\fesh v$-)finite
  element Hessian $\matfeop Hv$ to be the unique element of $\fesh h$
  that satisfies
  \begin{equation}\label{FEhess:def}
    \ltwop{\matfeop Hv}\varPhi_{\Wh}
    =
    \duality{\D^2v}\varPhi
    \Foreach
    \varPhi\in\fesh v,
  \end{equation}
  where the $\D^2$ is the generalised Hessian.  The finite element
  Hessian is thus the generalised Hessian's $\leb2({\Wh})$ representation in
  $\fesh v$. Notice that since $v\in\fesh v$, it's weak gradient is piecewise smooth, and thus has a boundary value, which, in particular, is also piecewise smooth; this results in the duality pairing on the right hand side of~(\ref{eq:def:generalised-Hessian}) being representable as a boundary integral.
\end{Def}
 \begin{Obs}[finite element Hessian for non finite element functions]
 We can define the ($\fesh v$-)finite element Hessian of a function
 $v\in H^2(\Wh)$, without any modification.
 \end{Obs} 
\begin{Obs}[symmetry of Hessians]
  For any $v\in\sobh1({\Wh})$ satisfying
  \begin{equation}
    \restriction{\grad
      v\Transpose{\normalto\W}}{\boundary\W}
    \in
    \psetmats{\dual{\sobh{1/2}(\boundary\W)}}dd,
  \end{equation}
  the generalised
  Hessian $\D^2v$ is symmetric, and so is the finite element Hessian
  $\matfeop Hv$.
\end{Obs}
\begin{Def}[finite element convexity after \citet{AguileraMorin:09:article:On-convex}]
A function $v\in \sobh 1(\Wh)$, such that its gradient's trace,
$\restriction{\grad v}{\boundary{\Wh}}$ is in
${\sobh{1/2}(\boundary{\Wh})'}$, is said to be strictly \indexemph{finite element
convex with respect to $\fesh v$}, concisely $\fesh v$-convex, if and only if
\begin{equation}
  \ltwop{\matfeop Hv}\varPhi_{\Wh}
  \in\spdmats[R]d
  \Foreach
  \varPhi\in\fesh v\take\setof0:\varPhi\geq0.
\end{equation}
Note that the test functions are such that they are nonnegative
everywhere and strictly positive on a set of positive measure.
\end{Def}
\subsection{Nonvariational finite element method (NVFEM) for the oblique derivative problem}
With these definitions in place it is possible to design a scheme
aimed at approximating $u$ satisfying problem~(\ref{genNDVobl}), by
seeking $(\fe u,\matfe H,c)\in\fesh v\times\fesh h\times\reals$ such that
\begin{equation}
  \label{nvfemalg}
  \begin{split}
    &
    \ltwop{
      \matfe H
    }{
      \varPhi
    }
    _{\Wh}
    +
    \inta{
      \grad \fe u
      (\grad \varPhi)\transposed 
    }
    _{\Wh}
    -
    \inta{
      \grad \fe u
      (\normalto{\Wh})\transposed 
      \varPhi
    }
    _{\boundary{\Wh}} 
    =
    \zeromat,
    \\
    &
    \ltwop{
      \mat A\frobinner\matfe H
      +
      {\vec b}\cdot\grad \fe u
      +
      c\fe u
    }
    \varPhi
    _{\Wh}
    +
    \ltwop{
      \vec\beta
      \inner
      \grad\fe u
    }
    \varPhi
    _{\boundary{\Wh}}+\ltwop{U}\kappa_{\Wh}+\ltwop{c}\varPhi_{\Wh}
    =
    \ltwop
    r\varPhi
    _{\Wh}
    +
    \ltwop
    s\varPhi
    _{\boundary{\Wh}}
  \end{split}
\end{equation}
for all $\varPhi\in\fesh v,\kappa\in\reals$. 

The nil sum constraint on $u$, the exact solution
of~(\ref{oblique:bvp}), needed to ensure its uniqueness is discretised by
seeking an additional unknown scalar (instead of directly including
this condition in the finite element space) $c$ as a Lagrange
multiplier, implemented by the inclusion of the following sum
\begin{equation}
  \label{eq:weak-Lagrange-multiplier-zero-sum-condition}
  \ltwop{\fe u}{\kappa}_{\Wh} 
  +
  \ltwop{c}{\fe \varPhi}_{\Wh} 
  =0
\end{equation}
in~(\ref{nvfemalg}).
Setting $\fe\varPhi = 0$ in~(\ref{nvfemalg}) gives us 
\begin{equation}
  \inta{\fe U}_{\Wh} = 0.
\end{equation}
Then, upon choosing $\varPhi\in\fesh v\meet\sobhz1(\Wh)$, we obtain
\begin{equation}
    \langle c,\varPhi\rangle_{\Wh}
    =
    \langle r-\mat A\frobinner\matfe H
    -
    \vec b\inner\grad\fe u-c\fe u,\varPhi\rangle_{\Wh}
\end{equation}
for all $\fe\varPhi\in\fesh v\meet\sobhz1(\Wh)$, which
tells us that $c$ is in fact the $\leb2({\Wh})$ projection
of
\begin{equation}
 r-\mat A\frobinner\matfe H
    -
    \vec b\inner\grad\fe u-c\fe u
\end{equation}
onto $\fesh v\meet\sobhz1(\Wh)$. Since $c$ is a constant, and the
only constant in $\fesh v\meet\sobhz1(\Wh)$ is zero, we deduce that
both integrals must be zero.

Note that the upper equation in~(\ref{nvfemalg}) is for a $(1,1)$ tensor on $\R d$, hence equivalent
to a system of $d^2$ equations, which, thanks to the symmetry of the
finite element Hessian, can be reduced to $d(d+1)/2$ equations; it is equivalent
to
\begin{equation}
  \matfe H=\matfeop H\fe U.
\end{equation}
The NVFEM, whose details for the Dirichlet boundary conditions are
described by \citet{LakkisPryer:11:article:A-finite}, can be viewed as
a mixed method, where we compute both the numerical solution $\fe u$
and its finite element Hessian $\matfe H=\matfeop H\fe U$, as an
auxiliary variable.  We stress, however that the variable $\matfe H$
becomes essential in nonlinear problems where the nonlinearity depends
on the Hessian. In fact, not only is accessing the finite element
Hessian necessary for the internal NVFEM algorithm, but as we see in
\S\ref{sec:nonlinear-solver}, it plays a crucial role in the outer
nonlinear solver and must therefore be returned by an implementation
of NVFEM.

Note also that (\ref{nvfemalg}) constitutes a departure from standard
FEMs in that the boundary condition is tested simultaneously with the
PDE, which is subsequently not integrated by parts.  It is therefore
not trivial that the solution of (\ref{nvfemalg}) should converge to the
exact solution of (\ref{genNDVobl}) in any meaningful sense.  We are
undertaking the analysis of this problem in a separate research.  The
numerical experiments we have conducted so far show that convergence
to optimal order can be obtained, at least for uniform meshes, if
the gradient recovery is used along side the Hessian recovery.
\section{A Newton--Raphson method for the Monge--Ampère with transport
  boundary condition}
\label{sec:nonlinear-solver}
In order to approximate $u$ satisfying to the nonlinear
problem~(\ref{eq:theprob}), we first work out the Newton--Raphson
method for the nonlinear problem, resulting in a sequence of solutions
$u_n$ to problems in the form of~(\ref{genNDVobl}) with $u$ replaced
by $u_n$.  As discovered by
\citet{LoeperRapetti:05:article:Numerical}, a Newton--Raphson
iteration, possibly with a damped stepsize converges to the exact
solution at the continuum level.  The main difficulty is to show that
the convexity of the Newton--Raphson iterate $u_n$ is preserved with
respect to $n$.  This leads to a sequence of well--posed elliptic
problems.  After discretisation with the finite element Hessian, it
turns out that the discrete problem inherits this property.
We now recap the results of~\citet{LakkisPryer:13:article:A-finite}
and then adapt them to problem~(\ref{eq:theprob}).
\subsection{Elliptic operators}
\label{def:elliptic-operators}
\label{ellipdef}
Consider a general Nemitsky-type (possibly nonlinear) operator of the form
\begin{equation}
  \begin{gathered}
    v\mapsto \cF[v]
    \\
    \cF[v(\vec x)]:=F(\vec x,v(\vec x),\grad v(\vec x),\D^2v(\vec x)),
  \end{gathered}
\end{equation}
which is well defined for functions $v\in\cont2(\W)$,
for some given (possibly nonlinear) function
\begin{equation}
  \funk F{\W\times\reals\times{\R d}\times\Symmatrices d}\reals
\end{equation}
where $\Symmatrices d$ indicates the vector space of symmetric linear
transformations on the Euclidean $\R d$.

Following \citet{CaffarelliCabre:95:book:Fully}, for an open set
$\mathscr{C}\subset\Symmatrices d$, the operator $\cF[\cdot]$
is called \indexemph{elliptic on $\cC$} if and only if, for each $(\xrp,\mat M)\in\W\times\R{d+1}\times\mathscr{C}$ there exist 
    $\lambda_\flat(\xrp,\mat M)\le\lambda_\sharp(\xrp,\mat M)$ in $\R+$, such that 
\begin{equation}
\begin{aligned}
  \label{elip}
  \lambda_\flat(\xrp,\mat M)
  \norm{\mat N}
  \le
  F(\xrp,\mat M+\mat N)-F(\xrp,\mat M)
  \le
  \lambda_\sharp(\xrp,\mat M)
  \norm{\mat N}
\end{aligned}
\end{equation}
for each $
  \mat N
  \in
  \Symmatrices d,$
where the matrix norm $\norm{\mat M}$ indicates the Euclidean-induced
operator norm (although the definition is independent of the choice of
norm except for the values of $\lambda_\flat$ and $\lambda_\sharp$).

If the largest possible set $\mathscr{C}$ for which~(\ref{elip}) is
satisfied is a proper subset of $\Symmatrices d$ we say that the
operator $\cF$ is \emph{conditionally elliptic}. The operator
$\cF[\cdot]$ is called \emph{uniformly elliptic on
  $\cC\subseteq\Symmatrices d$} if and only if
\begin{equation}
  \label{unifelip}
  0<
  \inf\nolimits_{%
    \W\times\R{1+d}\times\cC}
  \lambda_\flat%
  ,\tand
  \sup\nolimits_{%
    \W\times\R{1+d}\times\cC}
  \lambda_\sharp%
  <\infty
  ;
\end{equation}
the extremums defined by (\ref{unifelip}) are called \emph{lower} 
\index{lower uniform ellipticity constants}
and \indexemph{upper uniform
ellipticity constants}.  If the infimum in (\ref{unifelip}) is zero the operator
is called \indexemph{degenerate elliptic on $\cC$}.
\subsection{Smooth elliptic operators}
\label{def:smooth-elliptic}
If $F$ is differentiable~(\ref{unifelip}) can be obtained from
properties of the derivative of $F$. A generic $\mat M\in\realmats dd$
being written as
\begin{equation}
  \mat M
  =
  \squmatmd md
  ,
\end{equation}
the derivative of $F$ at $\mat M$ in the direction $\mat N$ is
represented by its $\grad_{\mat M}F(\xrp,\mat M)$, with respect to the Frobenius
product~(\ref{eq:def:Frobenius-inner-product}). Namely,
\begin{equation}\label{eq:def:Mderivative}
  \D_{\mat M}F(\xrp,\mat M)\mat N=:
  \grad_{\mat M}F(\xrp,\mat M)\frobinner\mat N
  \Foreach
  \mat N\in\realmats dd
\end{equation}
for some matrix $\grad_{\mat M}F(\xrp,\mat M)$, where
we have
\begin{equation}
  {\grad_{\mat M}F(\cdot,\mat M)
  =
  \providecommand{\klpfrac}[2]{\pd{\matentry m{#1}{#2}}F(\cdot,\mat M)}
  \dismatcommfromtofromto{\klpfrac}1d1d}.
\end{equation}
Usually, the function $F$ (and its gradient) are restricted to the
linear subspace $\Symmatrices d\subset\realmats dd$ in the $4$th
argument.  Therefore, if $F$ is differentiable then~(\ref{elip}) for
all $\mat{M}\in\cC$ is satisfied if and only if for each
$\mat{M}\in\mathscr{C}$ the matrix $\grad_{\mat M}F(\cdot,\mat M)$ is
(symmetric) positive definite, i.e.,
\begin{equation}
  \vec\xi\transposed\grad_{\mat M}F(\xrp,\mat M)\vec\xi
  \geq
  \lambda_\flat(\xrp,\mat M)\norm{\vec\xi}^2
  \Foreach \vec\xi\in\R d.
\end{equation}
Furthermore $\mathscr{C}=\Symmatrices d$ and $\lambda_\flat$ is independent of $\mat M$ 
if and only if the infimum condition in (\ref{unifelip}) is satisfied.
\begin{Lem}[ellipticity of the Monge--Ampère operator]
  \label{elliplemma} 
  The Monge--Ampère operator\footnote{
    Since the function $F$ generating the Monge--Ampère operator $\cF$
    does not depend on the values of the second variable representing the values
    of the operand ($v$ or $r$) we drop it.
  }
  \begin{equation}
    \label{eq:def:Monge-Ampere-operator}
    \cF[v]:=F(\vec x,\grad v,\D^2v)
    \text{ with }
    F(\xp,\mat M)
    :=
    \det\mat M
    -
    \frac{\rho(\vec x)}{\sigma(\vec p)}
  \end{equation}
  and $\rho$, $\sigma$ as described in \S\ref{sec:monge-ampere-problem}, is
  degenerate conditionally elliptic for $\mat M$ in the cone $\SPD(\R
  d)$ of symmetric positive definite linear transformations on $\R d$.
\end{Lem}
\begin{Proof}
  From the definitions in \ref{def:smooth-elliptic}, we need
  to show that $v\mapsto\det\D^2 v$ is elliptic.
  Recall the definition of the cofactor matrix, or tensor,
  of an invertible $\mat M$:
  \begin{equation}
    \label{eq:def:cofactor-matrix}
    \Cof\mat M
    :=
    \det(\mat M)
    \transinversemat M
    \eqncomment{  where $\transinversemat M:=\transposeof{\inversemat M}=\inverseof{\transposemat M}$}
  \end{equation}
  this definition can be extended by uniform continuity to singular matrices.
  By the definition of matrix invariants
  \citep{Bellman:97:book:Introduction} 
  we have, for each $\mat M,\mat N\in\realmats dd$ and $\theta\in\reals$,
  \begin{equation}
    \det\qp{\mat M+\theta\mat N}
    =
    \det\mat M
    +
    \Cof\mat M
    \frobinner
    \mat N
    \theta
    +
    \varrho(\theta)
  \end{equation}
  for a remainder function $\varrho$ satisfying
  \begin{equation}
    \abs{\varrho(\theta)}
    \leq
    \constdef[d]{const:Jacobi}
    \pownorm d{\mat M}
    \pownorm d{\mat N}
    \theta^2
    \Foreach\theta\in[0,1)
  \end{equation}
  for some $\constref[d]{const:Jacobi}$,
  from which we derive Jacobi's formula 
  \begin{equation}
    \label{eq:Cramers-rule}
    \D\det(\mat M)\mat N
    =
    \trace\qp{\Cof(\mat M)\mat N}
    =
    \Cof(\mat M)\frobinner\mat N
    \Foreach \mat M,\mat N\in\realmats dd.
  \end{equation}
  Thus, the gradient of $F$ with respect to the Frobenius
  inner product of matrices is
  \begin{equation}
    \grad_{\mat M}F(\xp,\mat M)
    =
    \Cof\mat M
    \Foreach 
    \mat M
    \in
    \realmats dd.
  \end{equation}
  This remains true when we restrict $F$ to matrices $\mat M$ (and
  variations thereof $\mat N$) in $\Symmatrices d$, or more
  specifically $\SPD(\R d)$.  Indeed, if $\mat{M}\in\SPD(\R d)$ then
  it is invertible, furthermore $\mat M^{-1}\in\SPD(\R d)$, and $\Cof\mat M=\det
  (M)\inversemat M$$\in\SPD(\R d)$. This holds because the eigenvalues of
  $\inversemat M$ are the reciprocals of the eigenvalues of $\mat M$,
  and since $\mat M$ is positive definite, all of its eigenvalues must
  be strictly positive. Thus for all $\vec\xi\in\R d$ we have that
  \begin{equation}
  \begin{split}
    \vec\xi\transposed\grad_{\mat M}F(\xp,\mat M)\vec \xi
    &
    =
    \det(\mat M)\vec \xi\transposed\inversemat M\vec \xi
    \\
    &
    \ge
    \frac
    {\pownorm2{\vec \xi}
    \det\mat M}{\lambda_\sharp},
  \end{split}
  \end{equation}
  where $\lambda_\sharp$ is the largest eigenvalue of $\mat M$.  Noting that
  since $\mat M$ is positive definite, its determinant is also
  strictly positive; it then follows that~(\ref{elip}) is
  satisfied. Since $\SPD(\R d)$ is a proper subset of $\Symmatrices d$
  this means that $\cF$ is only conditionally elliptic with maximal
  domain of ellipticity the functions whose Hessian is in $\SPD(\R
  d)$, i.e., the strictly convex functions. Finally noting that
  \begin{equation}
    \inf\nolimits_{\mat{M}\in\SPD(\R d)}\lambda_\flat(\mat{M})=0,
  \end{equation}
  it follows that $F$ is \indexemph{degenerate elliptic on} $\SPD(\R d)$.
\end{Proof}
\subsection{The Newton--Raphson method}
By Lemma~\ref{elliplemma} the operator $\cF[\cdot]$ is elliptic on
$\SPD(\R d)$. 
We introduce the cone of convex functions with nil sum on $\W$
\begin{equation}
  \cC:=
    \ensemble{v\in\cont2(\closure\W)}{\D^2 v(\vec x)\in\SPD(\R d)
    \Foreach \vec x\in\W
    \tand
    \inton v\W=0
    }
  .
\end{equation}
Furthermore, in order to capture the transport boundary
condition in (\ref{eq:theprob}), we introduce the nonlinear
operator 
\begin{equation}
  \label{eq:def:target-set-level-function}
  \mathscr{B}[u]:= b(\grad u)
  .
\end{equation}
With the notation from (\ref{eq:def:Monge-Ampere-operator}) and
(\ref{eq:def:target-set-level-function}), Problem~(\ref{eq:theprob})
consists in finding a function $\funk u\W\reals$ such that
\begin{equation}
  \label{mawopers}
  \begin{split}
    \cF[u(\vec x)] &= 0,\quad\vec x\in\W ,
    \\
    \mathscr{B}[u(\vec x)] &= 0,\quad\vec x\in\boundary\W .
  \end{split}
\end{equation}

To approximate the solution of~(\ref{mawopers}) we will apply the
\indexemph{Newton--Raphson method}. For each $n\in\NO$, assuming
$u_n\in\cC$ is given, the \indexemph{Newton--Raphson iteration}
consists in finding $u_{n+1}\in\cC$ satisfying
\begin{equation}
  \label{eq:exact-Netwon-Raphson-iterate}
  \begin{gathered}
    \D\cF[u_n(\vec x)](u_{n+1}(\vec x)
    -u_n(\vec x))+\cF[u_n(\vec x)] 
    = 0,
    \text{ for }\vec x\in\W ,
    \\
    \D\mathscr{B}[u_n(\vec x)](u_{n+1}(\vec x)-u_n(\vec x))
    +\mathscr{B}[u_n(\vec x)] 
    = 0,
    \text{ for }\vec x\in\boundary\W ,
  \end{gathered}
\end{equation}
where the $\D\cF$ and $\D\cB$ are the (infinite dimensional)
directional derivatives, explicitly calculated as
\begin{equation}
\begin{split}
  \D\cF[v]w  
  &:=\D F(\cdot,\grad v,\D^2v)(\vec 0,\grad w,\D^2w) \\
  &
  = 
  \Cof (\D^2v)\frobinner\D^2w
  +
  \frac{\rho}{\sigma(\grad v)^2}
  \D\sigma(\grad v)\grad w
  ,
\end{split}
\end{equation}
and
\begin{equation}
  \D\mathscr{B}[v]w
  :=
  \D b(\grad v)\grad w
  .
\end{equation}

It follows that at the $n$-th Newton--Raphson iteration we have to solve,
for the unknown $\theta_{n+1}:=u_{n+1}-u_n$, the oblique derivative elliptic
problem in nondivergence form (\ref{genNDVobl}) with the following
data
\begin{equation}
  \label{eq:def:linearised-Monge--Ampere-operators}
  \begin{aligned}
    \mat A(\vec x)
    &\gets
    \Cof\D^2u_n(\vec x)
    &&
    =:\hatmat A(\D^2u_n(x))
    ,
    \\
    \vec b(\vec x)
    &
    \gets
    \frac{\rho(\vec x)}{\sigma(\grad u_n(\vec x))^2}
    \grad\sigma(\grad u_n(\vec x))
    &&
    =:\hatvec b(\vec x, \grad u_n(\vec x))
    ,
    \\
    c(\vec x)&\gets0
    ,
    \\
    r(\vec x)
    &
    \gets
    -\det\D^2u_n(\vec x)
    +
    \frac{\rho(\vec x)}{\sigma(\grad u_n(\vec x))}
    &&
    =:\hat r(\vec x,\grad u_n(x),\D^2u_n(x))
    ,
    \\
    \vec\beta(\vec x)
    &\gets
    \grad b(\grad u_n(\vec x))
    &&
    =:
    \hatvec\beta(\grad u_n(\vec x))
    ,
    \\
    \tand
    &&&
    \\
    s(\vec x)
    &
    \gets
    -
    b(\grad u_n(\vec x))
    .
    &&
    =:
    \hat s(\grad u_n(x))
    .
  \end{aligned}
  \end{equation}
\section{The finite element scheme}
\label{sec:The scheme}%
We apply the NVFEM (\ref{nvfemalg}), to
approximate the terms $u_n$ of the sequence defined by
(\ref{eq:exact-Netwon-Raphson-iterate}).
\subsection{NVFEM--Newton--Raphson with plain finite element gradient}
A first attempt to discretise the Newton--Raphson iteration
(\ref{eq:exact-Netwon-Raphson-iterate}) can be derived,  as follows,
 for each $n\in\NO$, assuming $(\fe u_n,\mat
H_n)\in\fesh{v}\times\fesh h$ is given, find $(\fe u_{n+1},\mat
H_{n+1},c_{n+1})\in\fesh v\times\fesh h$ such that
\begin{equation}
  \begin{split}
    \label{suboptalg}
    \ltwop{\matfe H_{n+1}}{\fe\varPhi}_{\Wh}
    +
    \inta{\grad\fe u_{n+1}(\grad\fe\varPhi)\transposed}_{\Wh}
    +
    \ltwop{\grad\fe u_{n+1}{(\normalto{\Wh})\transposed}}{\fe\varPhi}_{\boundary{\Wh}}
    =
    \zeromat
    \\
    \Foreach\fe\varPhi\in\fesh v,
    \\
    \ltwop{\hatmat A(\matfe H_n)\frobinner(\matfe H_{n+1}-\matfe H_n)
      +\hatvec b(\cdot,\grad\fe u_n)\cdot\grad[\fe u_{n+1}-\fe u_n]
      +
      F(\cdot,\grad\fe u_n,\matfe H_n)}{\fe\varPhi}_{\Wh}
    \\
    +
    \ltwop{%
      \hatvec\beta(\grad\fe U_n)\cdot\grad[\fe u_{n+1}-\fe u_n]+\hat s(\grad \fe u_n(x))
    }{%
      \fe\varPhi}_{\boundary{\Wh}}
    +
    \ltwop{\fe u_{n+1}}{\kappa}_{\Wh}
    +
    \ltwop{ c_{n+1}}{\varPhi}_{\Wh}=0
    \\
    \Foreach\varPhi\in\fesh v,\kappa\in\reals.
  \end{split}
\end{equation}
\subsection{Shortcomings of the plain gradient approach of (\ref{suboptalg})}
Numerical experiments, show that algorithm
(\ref{suboptalg})--(\ref{eq:weak-Lagrange-multiplier-zero-sum-condition})
produces sequences that appear to be divergent for $\poly1$ elements.
Convergence is recuperated for $\poly k$ elements with $k\geq2$, but,
as the numerical experiments in Appendix~\ref{theappendix} show,
 convergence rates are suboptimal (in a function
approximation sense) in the $\leb 2({\Wh})$ norm. For instance, for
$\poly2$ elements, with the expected optimal convergence rate being
$3$, we observe a rate of $2$ at best.
\subsection{Boundary approximation}
We believe that the suboptimal results mentioned above caused by approximating a curved convex domain
by a polyhedral domain. The use of $\poly k$, $k\ge 2$
approximation requires the positioning of degrees of freedom on the
approximating boundary that in fact lie in the interior of the true
domain. This is why we observe a ``cap" on our convergence rates.  The
solution to this problem, at least from an empirical point of view,
based on extensive numerical computation is provided by the use of
\indexemph{gradient recovery}, in the case of $\poly 1$ elements (we still observe suboptimal rates in the $\leb 2({\Wh})$ norm for quadratics and higher).
\begin{Def}[projection-based gradient recovery] 
\par
We define the \indexemph{projection-based gradient recovery operator}
\begin{equation}
  \dfunkmapsto{\vecfeop g}v{\fesh v}{\vecfeop gv:=\vecfeop p\grad v}{\fesh g}
\end{equation}
where $\funk{\vecfeop p}{\setvecs{\leb 2({\Wh} )}d}{\setvecs{\fesh v}d}$ 
is the $\setvecs{\leb2({\Wh})}d$-projection operator. Explicitly this can be written as
\begin{equation}
  \label{eq:def:explicit-projection-based-gradient-recovery}
  \vecfeop gv\in\fesh g
  :
  \ltwop{
    \vecfeop gv
    -
    \grad v
  }{
    \fe\varPhi
  }
  _{\Wh}
  =0
  \Foreach
  \fe\varPhi\in\fesh v
  .
\end{equation}
Other gradient recovery operators, e.g., the one given by Zienkiewicz--Zhu,
which involves a more efficient local projection would be possible, but
we do not explore this issue in the current work.
\end{Def}
\subsection{FE Hessian with gradient recovery}
The standard FE Hessian operator, $\matfeop H$, defined
in~(\ref{FEhess:def}) is implemented in the NVFEM-Newton-Raphson by
it's inclusion in~(\ref{suboptalg}). Now that we are equipped with the
gradient recovery operator, $\vecfeop G$, given
by~(\ref{eq:def:explicit-projection-based-gradient-recovery}), we are
inclined to define a new finite element Hessian operator
$\tildematfeop H$, where one replaces the appearance of
$\grad\fe u$ in~(\ref{suboptalg}), with the recovered gradient
$\vecfeop G\fe u$, resulting in the following definition.
\begin{Def}[finite element Hessian with gradient recovery]
  We first define the gradient recovered generalised Hessian
  $\cH $, acting on $v\in\sobh1(\Wh)$ via
  \begin{equation}
    \langle \cH v|\varphi\rangle
    :=
    -\langle\vecfeop Gv\D\varphi\rangle_{\Wh}
    +
    \duality{\vecfeop Gv\normalto\W\transposed}\varphi_{%
      \dual{\sobh{1/2}(\boundary\Wh)\times\sobh{1/2}(\boundary\Wh)}}
    \Foreach\varphi\in\sobh1(\Wh).
  \end{equation}
Then, thanks to finite element conformity $\fesh v\subseteq\sobh1(\W)$,
we may define the finite element Hessian with gradient recovery
operator $\tildematfeop H$, acting upon $v\in\sobh1$ as follows
\begin{equation}
  \ltwop{
    \tildematfeop Hv
  }{\varPhi}_{\Wh}
  =
  \duality{\cH v}\varPhi
  \Foreach
  \fe\varPhi\in\fesh v.
\end{equation}
\end{Def}
\subsection{${\tildematfeop h}$ versus ${\matfeop h}$}
The use of the finite element Hessian with gradient recovery operator
is motivated by empirical observations that convergence properties are
superior for piecewise linear finite element approximation when using
$\tildematfeop h$ in conjuction with $\vecfeop g$, as opposed to
$\matfeop h$ with $\vecfeop g$.
\subsection{Gradient recovery for $\poly1$ elements}
Upon applying the gradient recovery operator $\vecfeop G$, defined
by~(\ref{eq:def:explicit-projection-based-gradient-recovery}), in
algorithm~(\ref{suboptalg}) for $\poly1$ element approximation we
observe that it does converge. Moreover, we observe optimal
convergence results in this case (see the first experiment in
Appendix~\ref{theappendix}).

The advantage of using piecewise linear polynomial approximation in
this case is that even if we approximate the curved convex domain with
a polyhedral domain, the degrees of freedom on the approximating
boundary in fact lie on the exact boundary, so in this case we would
expect to see optimal convergence rates.  This however, is no longer
possible for $\poly2$ elements and higher, as the boundary needs to be
approximated better to obtain full
convergence.
\subsection{Gradient recovery for $\poly k$, $k\geq2$}
The gradient of our approximate solution may be discontinuous (this
discontinuity can occur when the true solution lies \emph{outside} of
the finite element space), in discordance with that of the actual
solution, which is assumed to be continuous. To this end, we wish to
use a gradient recovery operator $\vecfeop g$, which has superconvergent
properties as noted by \citet{Zlamal:77:inproceedings:Some},
i.e., $\vecfeop g\fe u$ will converge faster to $\grad u$,
than the discrete gradient of our approximate solution $\fe u$. We
introduce the recovered gradient into our system as an auxiliary
variable to be solved for; as such, each component of $\vecfeop g\fe u$
will lie in the finite element space $\fesh v$.
\subsection{NVFEM--Newton--Raphson with finite element gradient recovery}
We incorporate the gradient recovery operator into our system, by
replacing $\grad\fe u_{n+1}$ with $\vecfeop g\fe u_{n+1}$
in~(\ref{suboptalg}).  This swap of roles in the discrete gradient
operator, implies a possible swap of the Hessian recovery operator
$\matfeop h$ with \indexemph{modified Hessian recovery operator}
$\funk{\tildematfeop h}{\fesh v}{\fesh h}$ for any $V\in\fesh v$,
\begin{equation}
    \ltwop{
      \tildematfeop H\fe v
    }{
      \fe\varPhi
    }
    _{\Wh}
    +
    \inta{
      \qp{\vecfeop g\fe v}
      \qp{\grad\fe\varPhi}\transposed 
    }
    _{\Wh}
    -
    \inta{
      \vecfeop g\fe v\,
      (\normalto{\Wh})\transposed 
      \fe\varPhi
    }
    _{\boundary{\Wh}} 
    =
    \zeromat
    \Foreach\fe\varPhi\in\fesh v.
\end{equation}
Rewriting the Newton--Raphson scheme \eqref{suboptalg} using $\tildematfeop H$
instead of $\matfeop H$, in \indexemph{incremental form} reads as follows,
for each $n\in\NO$,
\begin{enumerate}
\item 
  given \({
    (\fe u_n,\vecfe g_n,\matfe{H}_n)
    \in
    \fesh v
    \times
    \fesh g
    \times
    \fesh h
  }\), 
  satisfying
  \begin{equation}\label{NRGR:require}
    \begin{split}
      &\quad\vecfe g_n=\vecfeop g\fe u_n,
      \quad
      \matfe H_n={\tildematfeop h}\fe u_n,
      \\
      &\fe u_n\mbox{ is strictly finite element convex},
    \end{split}
  \end{equation}
\item 
  find the \indexemph{Newton--Raphson increment}~
  $\fe\varTheta\in\fesh{v}$ (along with its recovered gradient
  $\vecfeop G\fe\varTheta=:\vecfe{\varGamma}$ and its modified recovered Hessian
  $\tildematfeop H\fe\varTheta=:{\matfe\varDelta}$ and a scalar $c$) such that:
 \begin{equation}\label{NR:inc}
   \begin{split}
     \ltwop{\matfe\varDelta}{\fe\varPhi}_{\Wh}
     +
     \inta{\vecfe{\varGamma}\,\grad\fe\varPhi\transposed}_{\Wh}
     -
     \inta{\vecfe{\varGamma}\,{(\normalto{\Wh})\transposed}\fe\varPhi}_{\boundary{\Wh}}
     &=
     \zeromat
     \Foreach\varPhi\in\fesh v,
     \\
     \ltwop{\vecfe{\varGamma}}{\fe\varPhi}_{\Wh}
     -
     \ltwop{\grad\varTheta}{\fe\varPhi}\rangle_{\Wh}
     &=
     \vec0
     \Foreach\varPhi\in\fesh v,
     \\
     \ltwop{
       \hatmat A(\matfe H_n)
       \frobinner
       \matfe\varDelta
       +
       \hatvec b(\vecfe G_n)\inner\vecfe{\varGamma}
       +
       F(\cdot,\vecfe G_n,\matfe H_n)
     }{
       \fe\varPhi
     }_{\Wh}
     \\
     +
     \ltwop{
       \hatvec\beta(\vecfe G_n)
       \inner
       \vecfe{\varGamma}+\hat s(\vecfe G_n)}{\fe\varPhi}_{\boundary{\Wh}}
     +
     \ltwop{\varTheta}{\kappa}_{\Wh}
        +
        \ltwop{c}{\fe\varPhi}_{\Wh}
        &=
        0
        \Foreach\varPhi\in\fesh v,\kappa\in\reals.
   \end{split}
 \end{equation}
 where the functions $\hatmat A$, $\hatvec b$, $\hatvec\beta$, and
 $\hat s$ are given
 by~(\ref{eq:def:linearised-Monge--Ampere-operators}), with $\vecfe
 g_n$ and $\matfe{H}_n$ in place of $\nabla u_n$ and $\D^2u_n$,
 respectively,
  \item 
    define the \indexemph{next Newton--Raphson} iterate
    \begin{equation}\label{NRGR:iterate}
      \colvecthree{
        \fe u_{n+1}
      }{
        \vecfe g_{n+1}
      }{
        \qgroup{
          \matfe H
        }_{n+1}
    }
    :=
    \colvecthree{
      \fe\varTheta
    }{
      \vecfe \varGamma%
    }{
      \matfe \varDelta
    }
    +
    \colvecthree{
      \fe u_n
    }{
      \vecfe g_n
    }{
      \qgroup{\matfe H}_n
    }.
  \end{equation}
\end{enumerate}
\subsection{\fenics{} implementation}
We provide a pseudocode describing how we calculate the finite element
solution
of~\eqref{NRGR:require}--\eqref{NRGR:iterate}. The code
is implemented in \fenics, using a Newton--Raphson solver, where we
embed the first two linear
equations of~(\ref{NR:inc})
in the nonlinear map.  To do this, we first observe that although
the first linear equation in~(\ref{NR:inc}) is a tensor-valued equation
for \changefromto{$\widetilde{\matfe H}_\varTheta$}{\(\matfe\varDelta\)}
and $\vecfe \varGamma$, it can be collapsed
into the following equivalent scalar-valued equation
\begin{equation}
  \inton{
    \matfe\varDelta
    \frobinner
    \matfe\varXi
  }{\Wh}
  +
  \inton{
    \vecfe G_{\varTheta}
    \inner
    \divof{\matfe\varXi\transposed} 
  }{\Wh}
  -
  \inton{
    \vecfe G_{\varTheta}
    \inner
    (\matfe\varXi
    \normalto{\Wh})
  }{
    \boundary{\Wh}
  }
  =
  0
  \Foreach \matfe\varXi\in\fesh h,
\end{equation}
where the divergence of a matrix-valued map is taken row-wise (and produces a column):
\begin{equation}
  \getrow{\div\mat M}i=\sumifromto j1d\pd j\matentry mij
  \Foreach\qp{\text{row index}}i\integerbetween1d.
\end{equation}
Similarly, we collapse the vector-valued gradient recovery equation
into the equivalent scalar-valued equation
\begin{equation}
  \label{eq:alg:MAwithgradrec:grad:collapsed}
  \ltwop{
    \vecfe\varGamma
  }{
    \vecfe\varPsi
  }_{\Wh}
  -
  \ltwop{
    \grad\fe\varTheta
  }{
    \vecfe\varPsi
  }_{\Wh}
  =
  0
  \Foreach\vecfe\varPsi\in\fesh g
  .
\end{equation}
We may therefore include the linear components of~(\ref{NR:inc}) that
involve the gradient recovery, zero-average constraint and Hessian
recovery operations, in a \indexemph{global discrete nonlinear
  operator} $\funk{\feop{n}}{\fesh{Y}}{\fesh{Y}}$, where
\begin{equation}
  \label{Y:def}
  \fesh{Y}:=\fesh{v}\times\fesh{g}\times\fesh{h}\times\reals,
\end{equation}
implicity defined at a given
$\vecfour{\fe{u}}{\vecfe{g}}{\matfe{h}}c\in\fesh{Y}$, via the
$\leb2(\W)$-Riesz representation on $\fesh{Y}$, by
\begin{multline}
  \ltwop{
    \feop n
    (\fe u,\vecfe g,\matfe h,c)
  }{
    (\fe\varPhi,\vecfe\varPsi,\matfe\varXi,\kappa)
  }
  \\
  :=
  \inton{
    \qp{\vecfe g-\grad\fe u}\inner\vecfe\varPsi
  }
  {\Wh}
  +
  \inton{
    \matfe h\frobinner\matfe\varXi
  }{\Wh}
  +
  \inton{
    \vecfe g\inner(\divof{\matfe\varXi\transposed})
  }
 {\Wh} 
  +
  \inton{
    \vecfe g\inner(\matfe\varXi\normalto\W)
  }{
    \boundary{\Wh}}
  \\
  \ltwop{
    F(\cdot,\vecfe g,\matfe h)
  }{
    \fe\varPhi
  }
  \rangle_{\Wh} 
  +
  \ltwop{
    b(\vec Z)
  }{
    \fe\varPhi
  }
  _{\boundary{\Wh}}
  +
  \ltwop
  c{\fe\varPhi}
  _{\Wh} 
  +
  \ltwop{\fe u}\kappa
  _{\Wh} 
  \\
  \Foreach \qp{\fe\varPhi,\vecfe\varPsi,\matfe\varXi,\kappa}\in \fesh Y.
\end{multline}
In the Newton--Raphson method applied to 
solve $\feop n(\fe u,\vecfe g,\matfe h,c)=0$,
the $n$-th step reads as
follows:
\begin{equation}
  \label{alg:nonlinearMA}
  \begin{split}
    &
    \text{given }
    \vecfour{\fe{u}_n}{\vecfe{g}_n}{\matfe{h}_n}{c_n}\in\fesh{Y}, 
    \\
    &
    \text{find }
    \vecfour{\fe\varTheta}{\vecfe \varGamma}{\matfe \varDelta}{c_{n+1}}\in\fesh{Y}
    \text{ such that }
    \\
    &
    \ltwop{
      D\feop n
      \vecfour{
        \fe u_n}{
        \vecfe g_n}{
        \matfe h_n}{
        c_{n}
      }
      \vecfour{
        \fe\varTheta%
      }{%
        \vecfe\varGamma%
      }{%
        \matfe\varDelta%
      }{c_{n+1}}}{%
      \vecfour{\fe\varPhi}{\vecfe\varPsi}{\matfe\varXi}\kappa%
    }
    \\
    &
    \phantom-
    =    
    -\ltwop{%
      \feop n\vecfour{\fe{u}_n}{\vecfe{g}_n}{\matfe{h}_n}{c_n}}{%
      (\fe\varPhi,\vecfe\varPsi,\matfe\varXi,\kappa)%
    }
    \Foreach \qp{\fe\varPhi,\vecfe\varPsi,\matfe\varXi,\kappa}\in \fesh Y,
    \\
    &
    \text{define }
    \vecthree{\fe{u}_{n+1}}{\vecfe{g}_{n+1}}{\matfe{h}_{n+1}}%
    :=
    \vecthree{\fe{u}_n}{\vecfe{g}_n}{\matfe{h}_n}%
    +
    \vecthree{\fe\varTheta}{\vecfe\varGamma}{\matfe\varDelta}%
  \end{split}
\end{equation}
Note that $\feop n$, which depends on four finite-dimensional vectors,
is nonlinear only in the first variable while linear in the last three
variables. Hence the three equations corresponding to the three
derivatives in the ``linear variables'' are equivalent to gradient recovery,
zero-average constraint and Hessian recovery operations in~(\ref{NR:inc}),
whereas the equation corresponding to the first (nonlinear) variable
yields the Newton--Raphson linearisation of the nonlinear problem. The \fenics Newton--Raphson solver
that we used calculates the derivative of the nonlinear form, $\feop
n$, symbolically. It is however, possible to provide the solver with
the derivative, $\D\feop n$, manually, if needed.
\subsection{The linear system}
Each step of the Newton-Raphson method involves solving a linear
system (corresponding to a nonvariational linear elliptic equation
with and oblique derivative) of the form
\begin{equation}\label{lin:system}
  \nummat E
  \Transpose{
    \disrowvecfour{\transnumvec\theta}{\transnumvec\gamma}{\transnumvec\delta}{c}
  }
  =
  \numvec f,
\end{equation}
where $\nummat{E}$ is a square matrix of
$(1+\fracl32d+\fracl12d^2)N+1$ ($d$ is
the spatial dimension and $N=\dim\vespace$) in
and the vectors
\begin{equation}
  \numvec\theta\in\R N
  ,
  \numvec\gamma\in\R{dN}
  ,
  \numvec\delta\in\R{Nd(d+1)/2}
  c\in
  \vec{F}\in\R{1+(1+d+d(d+1)/2)N},
\end{equation}
quantify the finite element functions of the discrete Newton--Raphson increment
$(\fe\varTheta,\vecfe\varGamma,\matfe\varDelta)$ (and the Lagrange multiplier
$c_{n+1}$) appearing in (\ref{alg:nonlinearMA}). In particular,
\begin{equation}
  \label{theta:def}
  \varTheta(\vec x)
  :=
  \numvec\theta\transposed{\vecfe\varPhi}(\vec x)
  \text{ for }\vec x\in\W,
\end{equation}
where $\vecfe\varPhi=\seqidotsfromto{\vecfe\varPhi}1N$ denotes the (column)
vector of nodal basis functions of $\fesh v$. Similarly for the
(column of columns) $\numvec\gamma=\vecidotsfromto{\numvec\gamma}1d$
for the gradient's increment
$\vecfe\varGamma=\seqidotsfromto{\fe\varGamma}1d$, where each geometric
(physical) coordinate $\fe{\varGamma}_{\alpha}$,
$\alpha\integerbetween1d$, is associated with a vector
${\numvec\gamma}_\alpha\in\R N$, via
\begin{equation}\label{Gtheta:def}
  \irow{\fe\varGamma}\alpha(\vec x)
  :=\vecentry{\numvec\gamma}{\alpha}\transposed{\vecfe\varPhi}(\vec x)
  \text{ for each $\vec x\in\W$}.
\end{equation}
Similarly, with one more geometric index, for the (symmetric) Hessian's increment
\begin{equation}\label{Htheta:def}
  \matentry{\matfe\varDelta}\alpha\beta(\vec x)
  :=
  \transposeof{\numvecentry{\numvec\delta}{\alpha\,\beta}}
  \vecfe\varPhi(\vec x)
  \text{ for $\vec x\in\W$ and upper-triangular indexing }
  \alpha\integerbetween1d,\:
  \beta\integerbetween\alpha d.
\end{equation}
The final entry $c$ encodes the Lagrange multiplier $c_{n+1}$ corresponding
to the function's total mass from~(\ref{alg:nonlinearMA}).

Explicitly in the $2=d$ case we have
\begin{equation}
  (\numvec\theta,\numvec\gamma,\numvec\delta,c),
  =
  (\numvec\theta,\listitwo{\numvec\gamma},\listitwo{{\vecentry{\numvec\delta}1}},\vecentry{\numvec\delta}{1\,2},\vecentry{\numvec\delta}{2\,2},c)
\end{equation}
\begin{equation}
  \label{Arep}
  \nummat E = 
  \begin{bmatrix}
    \Diag\numvec d
    & {\nummat C}_1
    & {\nummat C}_2
    & {\nummat B}_{1\,1}
    & {\nummat B}_{1\,2}
    & {\nummat B}_{2\,2}
    & \numvec d
    \\
    \nummat A_1
    & \nummat M
    &\zeroentry
    &\zeroentry
    &\zeroentry
    &\zeroentry
    &\zeroentry
    \\
    \nummat A_2 & \zeroentry &\nummat M &\zeroentry &\zeroentry &\zeroentry &\zeroentry
    \\
    \zeroentry & {\nummat R}_1 & \zeroentry & \nummat M &\zeroentry &\zeroentry &\zeroentry
    \\
    \zeroentry &  {\nummat R}_2 & \zeroentry & \zeroentry &\nummat M &\zeroentry &\zeroentry
    \\
    \zeroentry & \zeroentry & {\nummat R}_2 & \zeroentry & \zeroentry & \nummat M & \zeroentry
    \\
    \transnumvec{d}
    &\zeroentry
    &\zeroentry
    &\zeroentry
    &\zeroentry
    &\zeroentry
    &\sum\numvec d  
  \end{bmatrix}
\end{equation}
with the blocks explicitly defined in pseudocode~(\ref{pseudoMA}),
and the right hand side $\numvec f$
\begin{equation}
  \label{RHS:def}
  \numvec f =
  {\vecidotsfromto{\numsca f}1N}\transposed = 
  -
  \ltwop{
    \feop n( \vecfour{\fe{u}_n}{\vecfe{g}_n}{\matfe{h}_n}{c_n})
  }{
    \vecfe\varPhi
  }_{\Wh}
  \tand
  \numvecentry fi=0
  \text{ for }i>N.
\end{equation}
The blocks appearing in (\ref{Arep}) are defined in \S\ref{pseudoMA},
which also summarises the whole procedure.
\begin{Alg}[Newton--Raphson--NVFEM-with-recovery]
  \label{pseudoMA}
  \begin{algorithmic}[1]
    \Require{%
      $\funk\rho\W\reals$
      $\funk\sigma\Y\reals$,
      $\funk b{\R d}\reals$,
      $\tol\in\R+$,
      $\itermax\in\naturals$,
      $\fesh v\gets\poly k(\mypart)\meet\cont0(\Wh)$
      Galerkin finite element space
      with basis $\vecfe\varPhi=\Transpose{\rowvecdotsfromto{\fe\varPhi}1N}$
      on $\mypart$ triangulation of $\W$,
      initial guess
      $\fe u_0\in\fesh v\times\reals$,
    }
    \Ensure{%
      $\fe U\in\fesh v$, $\vecfe g\in\fesh v^d$,
      $\vecfe h\in\Symmats[\fesh v]d$
      approximation of $u$ satisfying
      of (\ref{eq:monge-ampere-problem})
      , $\grad u$ and $\D^2 u$.
    }
    \Procedure{Newton--Raphson--NVFEM}{
      $\fe u_0,\rho,\sigma,b,\tol,\itermax$
    }
    \State{$n\gets 0$}\Comment{initialise the iteration counter}
    \State{$r\gets 1$}\Comment{initialise the Netwon--Raphson residual}
    \State{$c\gets 1$}\Comment{initialise the Lagrange multiplier}
    \State{\({%
      \nummat M
      \gets
      \langle
      \vec\varPhi{\vec\varPhi}\transposed
      \rangle_{\Wh}
      }\)}%
    \Comment{mass matrix}
    \State{\(
      \numvec d\gets \langle\vec\varPhi,1\rangle_{\Wh}
      \)}%
    \Comment{total mass}
    \State{%
      $\numvec u\gets\inverse{\nummat M}\ltwop{U_0}{\vecfe\varPhi}$
    }
    \Comment{initialise the potential DOF vector}
    \For{\(\alpha\integerbetween1d\)}
    \Comment{loop over the geometric directions}
    \State{
      \(
      \nummat A_\alpha
      \gets
      \inton{\vec\varPhi\transposeof{\partial_\alpha\vec\varPhi}}{\Wh}
      \)}%
    \Comment{discrete potential-to-$\alpha$th-derivative map}
    \State{\(
      \nummat{R}_\alpha
      \gets
      \nummat A_\alpha
      -
      \inton{
        \getrow{\normalto{\Wh}}\alpha{\vecfe\varPhi}
        \Transpose{\vecfe\varPhi}
      }{\boundary\Wh}
      \)}%
    \Comment{discrete potential-to-$\alpha$th-derivative map with boundary}
    \State{%
      ${\numvec g_\alpha}\gets\inverse{\nummat m}\nummat a_\alpha\numvec u$
    }
    \Comment{initialise the gradient DOF vectors}
    \For{\(\beta\integerbetween\alpha d\)}
    \State{$\numvec{h}_{\alpha,\beta}\gets\inverse{\nummat m}\nummat R_\alpha\numvec g_\beta$
    }
    \Comment{initialise the modified recovered Hessian DOF vectors}
    \EndFor
    \EndFor
    \For{\(\alpha\integerbetween1d\)}
    \State{%
      $\matfe{G}_\alpha\gets{\numvec g_\alpha\transposed\vecfe\varPhi}$
    }
    \Comment{initialise the gradient}
    \For{\(\beta\integerbetween\alpha d\)}
    \State{%
      $[\matfe{H}]^\beta_\alpha\gets{\numvec{h}_{\alpha,\beta}\transposed\vecfe\varPhi}$
    }
    \Comment{initialise the modified recovered Hessian}
    \EndFor
    \EndFor
    \While{$n\le\itermax$~\textbf{and}~$r>\tol$}
    \For{$\alpha = 1,\ldots,d$}
    \State{\({
        \nummat C_\alpha
        \gets
        \langle{
          \frac{\rho}{\sigma({\vecfe g})^2}
          \pd\alpha(\sigma(\vecfe g))\vec\varPhi
          \vec\varPhi\transposed
        }\rangle_{\Wh}
        +
        \langle
        \pd\alpha (b({\vecfe g}))\vec\varPhi{\vec\varPhi}\transposed
        \rangle_{\boundary{\Wh}}
      }\)}
    \For{$\beta = \alpha,\ldots,d$}
    \State{%
      \({%
        {\nummat B}_{\alpha\,\beta}
        \gets
        -\langle\getentryij{\Cof\matfe h}\alpha\beta
        \vec\varPhi{\vec\varPhi}\transposed\rangle_{\Wh}
      }\)}
    \EndFor
    \EndFor
    \State{construct $\nummat E$ given by~(\ref{Arep})}
    \State{construct $\numvec f$ given by~(\ref{RHS:def})}
    \State{%
      solve linear system
      \({%
        \displaystyle
        \nummat E
        \disrowvecfour{\numvec\theta}{\numvec\gamma}{\numvec\delta}c\transposed
        =\numvec f
      }
      \)%
    }
    \State{$\fe\varTheta\gets\vec\theta\transposed\vecfe\varPhi$}
    \Comment{Update the potential's increment}
    \For{\(\alpha\integerbetween1d\)}
    \State{%
      $\matfe{\varGamma}_\alpha\gets{\numvec \gamma_\alpha\transposed\vecfe\varPhi}$
    }
    \Comment{Update the gradient's increment}
    \For{\(\beta\integerbetween\alpha d\)}
    \State{%
      \(
        [\matfe{\varDelta}]^\beta_\alpha
        \gets
        \numvec{\delta}_{\alpha,\beta}\transposed\vecfe\varPhi
        \)}
    \Comment{Update the modified recovered Hessian's increment}
    \EndFor%
    \EndFor%
    \State{
      $(\fe U,\vecfe G,\matfe H)\gets\cdot+(\fe\varTheta,\vecfe\varGamma,\vecfe\varDelta)$
    }
    \Comment{Update solution by adding just computed increment}
    \State{$n\gets n+1$}
    \Comment{Update iterate counter}
    \State{$r\gets\Norm{\feop n(\fe u,\vecfe g,\matfe h,c)}_{\leb\infty({\Wh})}$}
    \Comment{Update Newton--Raphson residual}
    \EndWhile
    \State{\Return{$(\fe u,\vecfe g,\vecfe h)$}}
    \EndProcedure
  \end{algorithmic}
\end{Alg}
\section{Experiments}
\label{results}
In this section we report on the numerical experiments. Our freely
available code \citep{KaweckiLakkisPryer:18:url:GitHub} requires a
\fenics \citep{LoggMardalWells:12:book:Automated} installation.  In
each case, for data $\W$, $\Y$ (hence $b$), $f$ and $g$ corresponding
to a known \indexemph{benchmark solution}, $u$, of
(\ref{eq:monge-ampere-problem}) we compute a sequence of
approximations $\listidotsfromto{\fe U}1M$ on a sequence of meshes
$\listidotsfromto{\fepartition T}1M$, with corresponding meshsize
$h_m$ and finite element space $\fesh[m]v:=\poly
k({\mypart[m]})\meet\cont0(\Wh[m])$.

In these examples, the source domain, $\W$, coincides with the unit disk
in $\reals^2$, and the target domain, $\Y$, is either given in the first example by the
unit disk in $\reals^2$, and in the second example by the ellipse
\begin{equation}
  \setofsuch{(x,y)}{\frac{1}{4}x^2+\frac{1}{9}y^2\leq1}.
\end{equation}
For each fixed
experiment, to compute the sequence of \indexemph{experimental order
  of convergence}\index{$\EOC$} defined as
\begin{equation}
  \EOC_{m,\linspace X}:=\frac{\log(\norm{e_{m+1}}_{\linspace X}/\norm{e_m}_{\linspace X})}{\log(h_{m+1}/h_m)}
  \text{ for }m\integerbetween1{M-1}
\end{equation}
where $e_m:=\fe U_m-u$ is the error and $\linspace X$ a possible
seminorm among $\leb2(\W_{\fepartition T_m})$, $\sobh1(\W_{\fepartition T_m})$
or approximations thereof
where $\grad\fe u$, and $\D^2\fe u$ are respectively replaced by
$\vecfeop g\fe u$, and $\matfeop h\fe u$ or $\tildematfeop h\fe u$.
We empirically
observe \indexemph{optimal convergence rates} when implementing the $\poly1$
gradient recovery scheme~(\ref{NRGR:require})--(\ref{NRGR:iterate}),
that is our experimental results adhere to the following trends:
\begin{gather}
  \Norm{u-{\fe u_m}}_{\leb2{(\Wh[m])}}
  \le
  \constdef{l2-error}h_m^2,
  \\
  \norm{u-{\fe u_m}}_{\sobh1{(\Wh[m])}}
  \le
  \constdef{h1-error}h_m,
\end{gather}
for some $\constref{l2-error},\constref{h1-error}>0$ independent of $\fepartition T_m$.

In contrast, we observe \indexemph{suboptimal convergence} when
implementing either (\ref{suboptalg})--(\ref{eq:weak-Lagrange-multiplier-zero-sum-condition})
or (\ref{NRGR:require})--(\ref{NRGR:iterate}), when the polynomial
degree $k\ge2$, i.e., we observe the following:
\begin{gather}
  \Norm{u-{\fe u_m}}_{\leb2{(\Wh[m])}}
  \leq
  \constdef{const:suboptimal-poly2-20}
  h_m^2
  \\
  \norm{u-\fe u_m}_{\sobh1{(\Wh[m])}}
  \leq
  \constdef{const:suboptimal-poly2-21}
  h_m^2,
\end{gather}
in contrast to the optimal (\indexen{best approximation}) rates
\begin{gather}
  \Norm{u-{\fe u_m}}_{\leb2{(\Wh[m])}}
  \le
  \constdef{const:best-approximation-20}
  h_m^{k+1}
  \\
  \Norm{u-{\fe u_m}}_{\sobh1{(\Wh[m])}}
  \le
  \constdef{const:best-approximation-21}
  h_m^k
  ,
\end{gather}
where the latter are the convergence results one would expected for an
optimal numerical scheme. The most likely cause for the suboptimal
convergence is the piecewise linear approximation of domains with
curved boundary. This (non)variational crime is commented on, and
treated by the use of isoparametric finite elements,
in~\citet{Scott:73:phdthesis:Finite-element} (in the general context of finite element
approximation theory), and so we expect isoparametric elements to
overcome this problem.

Throughout our experiments, we also look at estimating the rates $r_1,r_2,r_3,r_4,\tilde r_4$ for
the following convergence estimates:
\begin{gather}
  \Norm{u-{\fe u_m}}_{\leb2{(\Wh[m])}}
  \le
  \constdef{l2-error:2}h_m^{r_1}
  \\
  \norm{u-{\fe u_m}}_{\sobh1{(\Wh[m])}}
  \le
  \constdef{h1-error:2}h_m^{r_2},
  \\
  \Norm{\grad u-{\vecfeop g\fe u_m}}_{\leb2(\Wh[m])}
  \leq
  \constdef{h1rec-error}h_m^{r_3},
  \\
  \Norm{\D^2 u-{\matfeop h\fe u_m}}_{\leb2(\Wh[m])}
  \leq
  \constdef{h2-error}h_m^{r_4},
  \\
  \Norm{\D^2 u-{\tildematfeop h\fe u_m}}_{\leb2{(\Wh[m])}}
  \leq
  \constdef{h2rec-error:2}h_m^{\tilde r_4},
\end{gather}
note that we numerically estimate the constants~\constref{l2-error:2},~\constref{h1-error:2},
in all experiments, the constants~\constref{h1rec-error}
and~\constref{h2rec-error:2} only when implementing~(\ref{NRGR:require})--(\ref{NRGR:iterate}),
and the constant~\constref{h2-error} otherwise.

Another characteristic worth mentioning is that of the recovered gradient's
superconvergence~\citep{ZhangNaga:05:article:A-new-finite}. When
implementing~(\ref{NRGR:require})--(\ref{NRGR:iterate}), our all of
our experiments
the recovered gradient outperforms the standard gradient; in some cases we even
observe that the recovered gradient error is consistently close to an entire order
higher than that of the standard gradient, e.g., in the $\poly1$
approximation.

The third series of numerical examples presented in
Appendix~\ref{theappendix}
are examples of image intensity transport
on one fixed uniform mesh. We transport (the negative of) a bitmap
image of Gaspard Monge, between two geometric objects. Namely,
the source domain, $\W$, is the unit square $(-1/2,1/2)^2$, which
corresponds to the ``space" that the original bitmap image of Monge occupies
and solve for the approximation of
problem~(\ref{eq:monge-ampere-domain-transport-condition})--(\ref{eq:monge-ampere-problem})
with the following density functions:
\begin{equation}
  \rho:=\left\{
  \begin{aligned}
    2 & \mbox{ if the pixel is white},\\
    1 & \mbox{ if the pixel is black}
  \end{aligned}
  \right.
\end{equation}
and the constant function
\begin{equation}
  \sigma\equiv \frac{1}{|\W|}\int_\W\rho.
\end{equation}
The resulting effect is for the white areas elements to be expanded
and the black ones to be compressed.  Reporting the transformation of
a uniform rectangular grid (not the computational grid) under the
gradient or the recovered gradient map renders the original bitmap
using recangles that are small in areas where the image is black and
large where the image is white. Note how the continuity of the
recovered gradient is useful in adding smoothness to the output grid.
The computational mesh is chosen to match the resolution of the
bitmap. Althought the function $\rho$ as defined here is
discontinuous, this is not an issue as there is only one mesh and we
only look at the possible use of MAOT solver as a way to encode image
information in a purely discrete fashion (hence the actualy $\rho$
could be continuous and we are just looking at a piecewise projection
of it).

\section{Conclusion}
\label{sec:conclusion}
We have presented a nonvariational finite element method for solving
the Monge--Ampère optimal transport problem. To our knowledge, while
the problem has been tackled with finite differences this is the first
with Galerkin type approximations.  The advantages of the Galerkin
approximation, over finite differences, is the ease of implementation
(we have just modified widely available packages, \fenics in our case,
but other ones may be used), the reasonable localisation of the method
(no need for wide stencils, e.g.) and a simple approximation at the
boundary.  Furthermore the use of finite elements allows for higher
order methods (which should be possible for isoparametric elements)
and, by using the gradient recovery, a continuous approximation of the
gradient of the solution, which is an excellent approximation for the
transport map $\grad u$ in the original Monge problem.

We empirically demonstrate the ease of implementation and robustness
of our method, as well as its ability to capture optimal error results
(in the $\poly1$ case), through a series of experiments. We also
provide an ``image processing'' example on how our method can be used
to construct monitor-function displaced girds.  This exhibits a step
forward in the area of mass transportation, and methods for both
linear and fully nonlinear elliptic equations with linear or nonlinear
oblique boundary conditions, as well as demonstrating the
applicability of variations of the nonvariational finite element
method introduced in \citet{LakkisPryer:13:article:A-finite,LakkisPryer:11:article:A-finite}.  The computational achievements of this paper are freely available
for reader's benefit on \cite{KaweckiLakkisPryer:18:url:GitHub}.

In terms of future research, the formulation of this method poses the
currently open question of existence and uniqueness of a solution to
the numerical scheme~(\ref{NRGR:require})--(\ref{NRGR:iterate}), as
the question of the derivation of optimal (or suboptimal) error
bounds. In order to achieve optimal error bounds for arbitrary
polynomial degree $k$, a potential avenue would be to incorporate the use
of isoparametric approximation of the computational domain.
\label{theappendix}
\begin{figure}
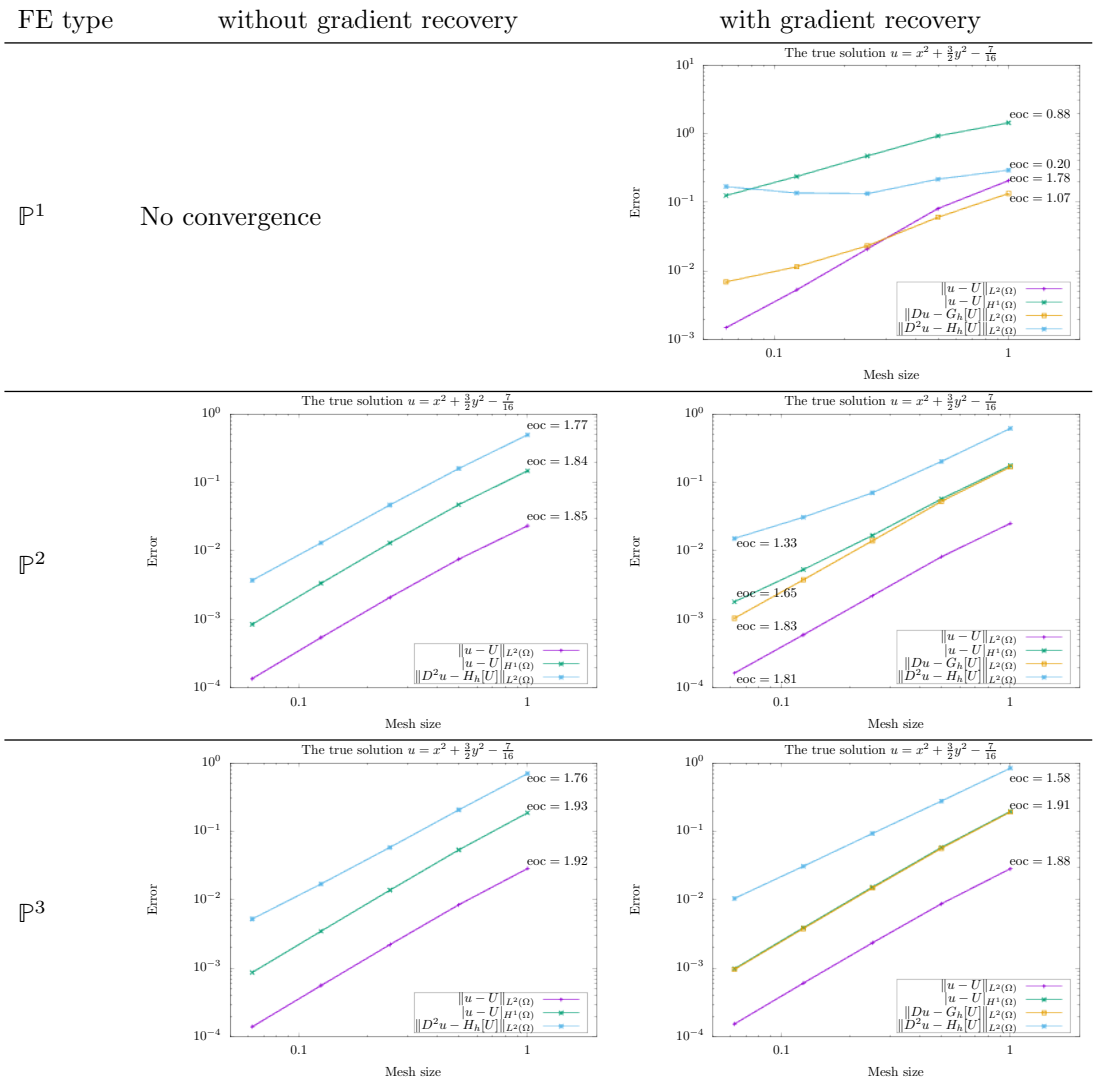

  \caption{In this case $\rho$ and $\sigma$ are chosen so that the true solution, $u(x,y) = x^2+3y^2/3-7/6$ and benchmark computations are performed without and with gradient recovery and with various polynomial degrees.}
  \providecommand{\benchbox}[2][0.5]{%
    \begin{minipage}{60mm}
      \scalebox{#1}{\input{Picture/#2}}
    \end{minipage}
  }
  \begin{center}
    \begin{tabular}{lcc}
      \toprule%
      FE type
      &
      without gradient recovery
      &
      with gradient recovery
      \\\midrule%
      $\poly1$
      &
    \begin{minipage}{60mm}
      \scalebox{1}{\input{Picture/noconvergence.tex}}
    \end{minipage}
      &
    \begin{minipage}{60mm}
      \scalebox{0.5}{\input{Picture/plot_exp_24.tex}}
    \end{minipage}
      \\\midrule%
      $\poly2$
      &
    \begin{minipage}{60mm}
      \scalebox{0.5}{\input{Picture/plot_exp_22.tex}}
    \end{minipage}
      &
    \begin{minipage}{60mm}
      \scalebox{0.5}{\input{Picture/plot_exp_25.tex}}
    \end{minipage}
      \\\midrule%
      $\poly3$
      &
    \begin{minipage}{60mm}
      \scalebox{0.5}{\input{Picture/plot_exp_23.tex}}
    \end{minipage}
      &
    \begin{minipage}{60mm}
      \scalebox{0.5}{\input{Picture/plot_exp_26.tex}}
    \end{minipage}
      \\\bottomrule%
    \end{tabular}
  \end{center}
\end{figure}
\clearpage
\providecommand{\mongeportrait}[3][0.88]{%
  \begin{minipage}{.45\textwidth}
    \begin{center}
      \includegraphics[width=#1\textwidth,trim=60 200 60 200,clip]{%
        MA-fenics-experiments/All-experiments/imagetransportexperiments/exp_#2/exp_#2_post.pdf}%
      \\
      #3
    \end{center}
  \end{minipage}
}
\begin{figure}
  \begin{center}
    \caption{
      \label{fig:monge}
      Gaspard Monge's mesh-portrait obtained by mass transporting a
      uniform rectangular mesh into a ``monitor'' function.}
    \begin{tabular}{rr}
      \begin{minipage}{.45\textwidth}
        \begin{center}
          \includegraphics[width=0.72\textwidth,trim=200 300 200 250,clip]{%
            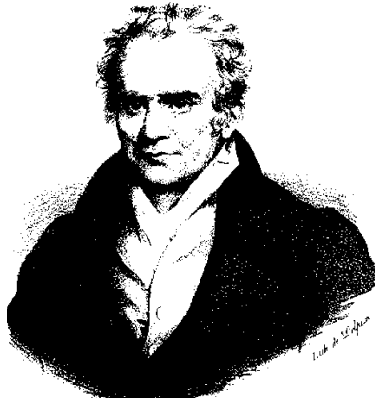}
          \\
          A bitmap of a portrait of Gaspard Monge,
          Lithography by F.S.~Delpech (Public Domain)
        \end{center}
      \end{minipage}
      &
      \mongeportrait[1]{17}{$\poly2$ FE without gradient recovery}%
      \\
      \mongeportrait{20a}{$\poly1$ FE with gradient recovery}
      &
      \mongeportrait{20b}{$\poly2$ FE with gradient recovery}
      \\
      \mongeportrait{18}{$\poly2$ FE without gradient recovery}%
      &
      \mongeportrait{19}{$\poly2$ FE with gradient recovery}%
    \end{tabular}
  \end{center}
\end{figure}
\clearpage
\bibliographystyle{abbrvnat}
\printindex%
\end{document}